\newif\iffinal
\finaltrue	

\documentclass[letterpaper,11pt,reqno]{amsart}
\RequirePackage[utf8]{inputenc}

\usepackage[usenames,dvipsnames,svgnames,table]{xcolor}

\usepackage[portrait,margin=3cm]{geometry}
\iffinal\else\usepackage[notref,notcite]{showkeys}\fi
\usepackage{mathrsfs,soul}
\usepackage{hyperref}
\usepackage[foot]{amsaddr}
\usepackage{amssymb,amsthm,amsfonts,amsbsy,latexsym,dsfont}
\usepackage[textsize=small]{todonotes}       
\usepackage{graphicx}
\usepackage[numeric,initials,nobysame]{amsrefs}
\usepackage{upref,setspace}
\usepackage{caption}
\usepackage{mathtools}

\usepackage{enumerate}
\newenvironment{enumeratei}{\begin{enumerate}[\upshape (i)]}{\end{enumerate}}
\newenvironment{enumeratea}{\begin{enumerate}[\upshape (a)]}{\end{enumerate}}

\newenvironment{enumerateA}{\begin{enumerate}[\upshape (A)]}{\end{enumerate}}

\numberwithin{equation}{section}
\numberwithin{figure}{section}
\numberwithin{table}{section}

\newtheorem{thm}{Theorem}[section]
\newtheorem{lem}[thm]{Lemma}
\newtheorem{cor}[thm]{Corollary}
\newtheorem{prop}[thm]{Proposition}
\newtheorem{defn}[thm]{Definition}

\newtheorem*{ass*}{Assumption}

\newtheorem{lemma}[thm]{Lemma}

\theoremstyle{definition}
\newtheorem{rem}{Remark}

\renewcommand{\leq}{\leqslant}
\renewcommand{\geq}{\geqslant}

\newcommand{\ind}{\mathds{1}}
\newcommand{\eps}{\varepsilon}

\newcommand{\set}[1]{\left\{#1\right\}}

\newcommand{\probc}{\stackrel{\mathrm{P}}{\longrightarrow}}

\newcommand{\weakc}{\stackrel{\mathrm{w}}{\longrightarrow}}
\newcommand{\convas}{\stackrel{\mathrm{a.e.}}{\longrightarrow}}

\def\qed{ \hfill $\blacksquare$}

\newcommand{\cB}{\mathcal{B}}
\newcommand{\cD}{\mathcal{D}}\newcommand{\cF}{\mathcal{F}}
\newcommand{\cG}{\mathcal{G}}\newcommand{\cI}{\mathcal{I}}
\newcommand{\cL}{\mathcal{L}}
\newcommand{\cM}{\mathcal{M}}
\newcommand{\cP}{\mathcal{P}}
\newcommand{\cT}{\mathcal{T}}

\newcommand{\vp}{\mathbf{p}}

\newcommand{\vw}{\mathbf{w}}

\newcommand{\mvtheta}{\boldsymbol{\theta}}

\newcommand{\mvpi}{\boldsymbol{\pi}}

\newcommand{\bL}{\mathbb{L}}
\newcommand{\bM}{\mathbb{M}}
\newcommand{\bR}{\mathbb{R}}

\newcommand{\sD}{\mathscr{D}}

\newcommand{\sG}{\mathscr{G}}

\newcommand{\sR}{\mathscr{R}}

\newcommand{\sZ}{\mathscr{Z}}

\DeclareMathOperator{\E}{\mathds{E}}
\DeclareMathOperator{\pr}{\mathds{P}}

\DeclareMathOperator{\var}{Var}

\DeclareMathOperator{\dist}{dist}

\DeclareMathOperator{\pre}{pre}
\DeclareMathOperator{\post}{post}

\newcommand{\sss}{\scriptscriptstyle}

\newcommand{\age}{{\sf Age} }
\newcommand{\epoch}{{\sf Epoch} }

\newcommand{\timea}{{\sf TimeAlive} }

\newcommand{\ac}{{\sf AC} }
\newcommand{\bc}{{\sf BC} }

\newcommand{\BP}{{\sf BP} }

\newcommand{\bcpmf}{_t h}
\newcommand{\acpmf}{h_t}

\newcommand{\out}{{\sf out} }
\newcommand{\good}{{\sf good} }
\newcommand{\bad}{{\sf bad} }
\newcommand{\birth}{{\sf Birth} }
\newcommand{\stt}{{\sf st} }

\begin{document}

\title[Change point detection in Networks]{Change point detection in Network models:\\ Preferential attachment and long range dependence}

\date{}

\author[Bhamidi]{Shankar Bhamidi}
\address{Department of Statistics and Operations Research, 304 Hanes Hall, University of North Carolina, Chapel Hill, NC 27599}
\author[Jin]{Jimmy Jin}
\author[Nobel]{Andrew Nobel}
\email{bhamidi, jimmyjin, nobel,  @email.unc.edu }

\maketitle
\begin{abstract}
Inspired by empirical data on real world complex networks, the last few years have seen an explosion in proposed generative models to understand and explain observed properties of real world networks, including  power law degree distribution and ``small world'' distance scaling. In this context, a natural question is the phenomenon of {\it change point}, understanding how abrupt changes in parameters driving the network model change structural properties of the network. We study this phenomenon in one popular class of dynamically evolving networks:  preferential attachment models. We derive asymptotic properties of various functionals of the network including the degree distribution as well as maximal degree asymptotics, in essence showing that the change point does effect the degree distribution but does {\bf not} change the degree exponent. This provides further evidence for long range dependence and sensitive dependence of the evolution of the process on the initial evolution of the process in such self-reinforced systems.  We then propose an estimator for the change point and prove consistency properties of this estimator. The methodology developed  highlights the effect of the non-ergodic nature of the evolution of the network on classical change point estimators.

\end{abstract}

\section{Introduction}
\label{sec:intro}

Motivated by the availability of data on many real world systems, the last few years have witnessed an explosion in both methodological as well as theoretical development of various complex network models. The aim of these models is to explain structural features observed in the data (e.g. power law degree distribution or ``small world'' connectivity) as well as understand and predict the behavior of dynamic processes on these networks including disease contact networks, search algorithms, random walks, evolution and dissolution of communities and a wide array of related processes  \cite{bollobas2001random,durrett-rg-book,van2009random,chung2006complex,newman2003structure,newman2010networks,albert2002statistical,dorogovtsev2002evolution}.  One sub-field of this vast field which has been particularly active is temporal or time varying networks. See the recent surveys \cite{boccaletti2014structure,holme2012temporal} and the references therein for both methodological developments as well as applications in a wide array of fields ranging from social networks and online communication, cell biology including temporal properties of protein interaction networks, and infrastructure systems such as the power grid.

Many such proposed models are driven by a collection of parameters that describe the evolution of the network.  A natural question in this context is the study of \emph{change points}, the effect of abrupt changes in parameters driving the evolution of the network, on structural properties of the network. To fix ideas, first consider the simplest version of the classical (offline) change point detection in the context of \emph{iid} data described as follows. Fix two distribution functions $F$ and $G$ (unknown but different) and a parameter $\gamma \in (0,1)$. Consider a stream of data $\set{X_i: 1\leq i\leq n}$ with distribution: for $i\leq \lfloor n\gamma \rfloor$ $X_i$ are \emph{iid} with distribution $F$ whilst for $i> \lfloor n\gamma \rfloor$, $X_i$ are \emph{iid} with distribution $G$ (and independent of the initial segment). Based on the observed data, $\set{X_i:1\leq i\leq n}$, the aim is then to estimate the change point $\gamma$ using estimators that are consistent as the sample size $n \to \infty$. 

In this spirit, this paper has two main goals:

\begin{enumeratea}
	\item We start with a variation of the standard preferential attachment model of evolving networks that incorporates a change point. This conceptually simple model allows for a simple interpretation of the effect of the change point on network dynamics. We rigorously study the effect of this change point on structural properties of the network including the scale-free or heavy tailed nature of the limiting degree distribution as well as asymptotics for the maximal degrees. 
	\item We then propose and study consistency properties of offline estimation procedures to detect the location of this change point from observed data. In particular this allows one to gain insight into the effect of the non-stationary nature of the evolution of the network model on various known heuristics for estimation in the \emph{iid} setting. 
\end{enumeratea}

\subsection{Organization of the paper}
\label{sec:org}
Both change point detection as well as preferential attachment models have witnessed enormous amount of work over the last few decades. We defer a fuller discussion of these two fields, their relevance to this paper as well as related work to Section \ref{sec:disc}. We start in Section \ref{sec:model} by defining the model. In Section \ref{sec:prelim} we setup notation required for the main results. Section \ref{sec:res} contains our main results,  starting with Section \ref{sec:res-deg-dist} that describes asymptotics for functionals of the networks including the degree distribution as well as maximal degrees as the network size $n\to\infty$.  Section \ref{sec:res-stat} formulates estimators to find the change point and formulates consistency properties for these estimators. Proofs for asymptotics of network functionals can be found in Section \ref{sec:proofs}. Section \ref{sec:leafclt-proofs} develops a functional central limit theorem for a specific functional of the network. Section \ref{sec:consistency} then uses this CLT to prove consistency of the proposed estimator. 

\subsection{Model formulation}
\label{sec:model}

We start by describing the original model of preferential attachment with no change point \cite{barabasi1999emergence,yule1925mathematical,simon1955class}. There are many variants of this model.  Throughout the paper we will consider the simplest case where the network at each stage is a tree. The methodology can be generalized to the general network setup.  Start with a single vertex at time $m=1$ (this vertex will be referred to as the \emph{root} or the original progenitor of the process and denoted by $\rho$). Fix a parameter $\alpha \geq 0$.  At each discrete time point  $1 < m\leq n$ a new vertex enters the system with a \emph{single edge} which it will then connect to a pre-existing vertex.  The vertex connects to a pre-existing vertex $v$ with probability proportional to the current degree of $v$ $+\alpha$. Let $\cT_m$ denote the graph at time $m$ and $\set{\cT_m: 1\leq m\leq n}$ be the entire graph valued process. Note that since each new vertex has one edge which it uses to connect to the current graph, $\cT_m$ for any $m$ is a tree (which we view as rooted at $\rho$). Thus for $m> 1$, the degree of every vertex is at least $1$. Further calling the pre-existing vertex that a new vertex attaches to as the \emph{parent} of this vertex, one can view this process as generating a directed tree with edges pointed from parents to children. 

We will soon switch over to a continuous time version of the process where it is convenient to work with a slight variant of the above process. Note that for a (directed rooted) tree, the degree of every vertex other than the root is $1+$ out-degree of the vertex. For the root, the degree and the out-degree coincide.  Now fix a single vertex at time $m=1$ and a parameter $\alpha >0$. The variant considered in this paper is as follows: at each stage $m> 1$ a new vertex enters the system and connects to a pre-existing vertex $v\in \cT_{m-1}$ with probability proportional to $1+\alpha+$ out-degree of $v$ in $\cT_{m-1}$. This variant results in all the same asymptotic properties as the original model and is slightly easier to deal with rigorously. 

This model has been studied extensively and in particular it is known \cite{Bollobas:2001:DSS:379831.379835} that the degree distribution converges in the large network limit. Precisely, for fixed $k\geq 1$, let $N_n(k)$ denote the number of vertices with degree $k$ in $\cT_n$. Then,
\begin{equation}
\label{eqn:pk-alpha-def}
\frac{N_n(k)}{n} \convas p_{\alpha}(k), \qquad \mbox{ where }	p_\alpha(k):= (2+\alpha)\frac{\prod_{j=1}^{k-1}(j+\alpha)}{\prod_{j=3}^{k+2}(j+2\alpha)}. 
\end{equation}
Here for $k=1$, we use the notation $\prod_{j=1}^{k-1} = 1$. Write $D_\alpha$ for a random variable with the above distribution. It is easy to check that there exists a constant $c >0$ such that 
\begin{equation}
\label{eqn:tail-bound}
	\pr(D_\alpha \geq k) \sim \frac{c}{k^{\alpha+3}}, \qquad \mbox{ as } k\to \infty. 
\end{equation}

 Further, arranging the degrees in $\cT_n$ in decreasing order as $M_n(1) \geq M_n(2) \geq \cdots M_n(n)$, it is known \cite{mori2007degree,bhamidi2012spectra} that for any fixed $k\geq 1$, there exists a non-degenerate probability distribution $\nu_k^\alpha$ on $\bR_+^k$ such that 
\begin{equation}
\label{eqn:max-degree-conv}
\left(\frac{M_n(j)}{n^{\frac{(1+\alpha)}{(2+\alpha)}}}: 1\leq j\leq k\right) \weakc \nu_\alpha^k.
\end{equation}

\subsubsection{\bf Model with change point:} 
\label{sec:model-cp-def}
Now fix two attachment parameters $\alpha, \beta > 0$, a change point parameter $\gamma \in (0,1)$, and a system size $n> 1$. The model does preferential attachment as before, but now the attachment dynamics changes after time $\lfloor n\gamma \rfloor$ namely 
\begin{enumeratea}
	\item For time $0 < m\leq \lfloor n\gamma \rfloor$, the new vertex entering the system at time $m$ connects to pre-existing vertices with probability proportional to their current out-degree $+1+\alpha$. 
	\item For time $\lfloor n\gamma \rfloor < t\leq n$, the new vertex connects to pre-existing vertices with probability proportional to their current out-degree $+1+\beta$. 
\end{enumeratea}
Let $\mvtheta=(\alpha,\beta,\gamma)$ be the driving set of parameters of the model.  We will let $\cT_{\mvtheta, m}$ denote the rooted tree at time $m$ and $\set{\cT_{\mvtheta,m}: 1\leq m\leq n}$ for the entire graph valued process. When the context is clear, for ease of notation we suppress the dependence on $\mvtheta $ and write  $\set{\cT_{m}: 1\leq m\leq n}$. This model is the main object of interest for the rest of the paper.   

\subsection{Preliminary notation}
\label{sec:prelim}
To state our main results we will need to define some additional objects. Recall the parameter set $\mvtheta := (\alpha, \beta, \gamma)$ used to construct the model. Let $\set{E_\alpha(k):k\geq 1}$ be a sequence of independent exponential random variables such that for each fixed $k\geq 1$, $E_\alpha(k)$ has rate $k+\alpha$. View the above as the inter-arrival times  of a point process $\cP_{\alpha}$ on $\bR_+$. More precisely write, 
\[L_\alpha(m) = E_\alpha(1)+\cdots + E_\alpha(m), \qquad m\geq 1.\]
Consider the point process 
\begin{equation}
\label{eqn:cp-alpha-def}
	\cP_\alpha:= (L_\alpha(1), L_\alpha(2),\ldots).
\end{equation}
Analogously define $\set{E_\beta(k):k\geq 1}$, $\set{L_\beta(k):k\geq 1}$ and the corresponding point process $\cP_\beta$. For fixed $t\geq 0$, write $N_\alpha(t):=\cP_{\alpha}[0,t]$ for the number of points in $\cP_\alpha$ which fall in the interval $[0,t]$. 

We will need variants of the above point process. Fix $j\geq 1$ and $\alpha >0$. Let $\cP_\alpha^j$ be the point process where we use the sequence of points $\set{E_\alpha(m): m\geq j}$ to construct the point process so that the first point arrives after an exponential rate $j+\alpha$ amount of time, the second point arrives at rate $j+1+\alpha$ after the first point and so forth. As before let $N_\alpha^j(\cdot)$ be the corresponding counting process and note that $N_\alpha^1(\cdot) = N_\alpha(\cdot)$. 

Define the constant
\begin{equation}
\label{eqn:a-def}
	a=\frac{1}{2+\beta}\log{\frac{1}{\gamma}}.
\end{equation}  
On the interval $[0,a]$, define the ``truncated'' exponential distribution described via the cumulative distribution function 
\begin{equation}
\label{eqn:ga-cdf-def}
	G_a(s) = \frac{1-\exp(-(2+\beta)s)}{1-\exp(-(2+\beta)a)}, \qquad s\in [0,a].
\end{equation} 
Write $\age$ for a random variable with distribution $G_a$ (the reason for this terminology will become clear in the proof). Generate a counting process $N_\beta(\cdot)$  as above (independent of $\age$) and let $X_{\ac} = N_\beta[0,\age]$, namely the number of points that occur before the random time $\age$. Here $\ac$ is a mnemonic for ``after change point''. 

We are now in a position to define the limiting degree distribution.   Consider the following integer valued random variable $D_{\mvtheta}$:
\begin{enumeratea}
	\item With probability $1-\gamma$, $D_{\mvtheta} = 1+X_{\ac}$. 
	\item With probability $\gamma$, $D_{\mvtheta}  = D_\alpha + N_\beta^{D_\alpha}[0,a]$ where $D_\alpha$ is a random variable with distribution as in \eqref{eqn:pk-alpha-def}, namely the limiting degree distribution {\it without} change point. More precisely, generate $D_\alpha$ with distribution as in \eqref{eqn:pk-alpha-def}. Conditional on $D_\alpha$, generate the point process $N_\beta^{D_\alpha}$ and count the number of points in the interval $[0,a]$ and add this to the original random variable $D_\alpha$. 
\end{enumeratea}
Write $\vp_{\mvtheta}=(p_{\mvtheta}(k): k\geq 1)$ for the probability mass function of the above random variable namely 
\begin{equation}
\label{eqn:p-theta-def}
	p_{\mvtheta}(k) = \pr(D_{\mvtheta} = k), \qquad k\geq 1.
\end{equation}

\section{Results}
\label{sec:res}

Let us now describe our main results. We state results about the asymptotic degree distribution in Section \ref{sec:res-deg-dist}. We formulate statistical procedures to estimate the change point and the associated consistency results in  Section \ref{sec:res-stat}. 

\subsection{Asymptotics for the degree distribution}
\label{sec:res-deg-dist}
Fix $\mvtheta \in \bR_+^2 \times (0,1)$. For fixed $k\geq 1$ let $N_n(k)$ denote the number of vertices with degree $k$ in the random tree $\cT_n$ constructed in the change point model as in Section \ref{sec:model-cp-def}. The random variable $D_{\mvtheta}$ in the following result is as defined in \eqref{eqn:p-theta-def}.  

\begin{thm}\label{thm:deg-dist}
 Fix $k\geq 1$. As $n\to\infty$ the degree distribution satisfies, 
 \[\frac{N_n(k)}{n} \probc \pr(D_{\mvtheta} = k), \]
 Further for $\alpha \neq \beta$ and $\gamma\in (0,1)$, $\vp_{\mvtheta} \neq \vp_\alpha$.  However there exist constants $0 < c < c^\prime $ such that for all $k\geq 1$
 \begin{equation}
 \label{eqn:mvtheta-tail}
 	\frac{c}{k^{\alpha+2}} \leq \pr(D_{\mvtheta} \geq k) \leq  \frac{c^\prime}{k^{\alpha+2}}.
 \end{equation}
\end{thm}

\begin{rem}
	This theorem says that one {\bf does feel} the effect of the change point in the empirical degree distribution if $\alpha \neq \beta$ and $\gamma \in (0,1)$, however comparing \eqref{eqn:mvtheta-tail} with \eqref{eqn:tail-bound}, for any fixed $\gamma \in (0,1)$, this does {\bf not} change the tail behavior. This is a little surprising as one might assume,  especially for $\gamma$ close to zero and $\beta < \alpha$ (where the no change point dynamics with $\beta$ instead of $\alpha$ results in a degree distribution with a heavier tail), the tail of the degree distribution might scale like $k^{-(3+\beta)}$, namely the dynamics of attachment driven by $\beta$ should kick in.  However this is not the case.   
\end{rem}

\begin{rem}
	The techniques developed in this paper easily extend to the setting of multiple change points. We describe these extensions in Theorem \ref{thm:multi-change-pt}. 
\end{rem}

The next result deals with maximal degree asymptotics. As before arrange the degrees in $\cT_n$ in decreasing order as $M_n(1)\geq M_n(2)\geq \cdots M_n(n)$.

\begin{thm}
\label{thm:max-deg}
 Fix $k\geq 1$ and consider the $k$ maximal degrees $(M_n(j):1\leq j\leq k)$. Then the sequence of $\bR_+^k$ valued random variables defined by setting 
\[\bM_n(k):= \left(\frac{M_n(j)}{n^{\frac{(1+\alpha)}{(2+\alpha)}}}: 1\leq j\leq k\right),\qquad n\geq 1, \]
is tight and bounded away from zero. 	
\end{thm}

\begin{rem}
 Comparing the scaling of the maximal degrees above to the setting of no change point as described in \eqref{eqn:max-degree-conv}, one sees that the maximal degrees do not feel the effect of the change point, at least in terms of their order of magnitude.  We further conjecture that $\set{\bM_n(k):n\geq 1}$ converge weakly to a non-degenerate distribution on $\bR^k_+$. We have not pursued this further in this paper. 
\end{rem}

\begin{figure}[htbp]
	\centering
		\includegraphics[scale=1]{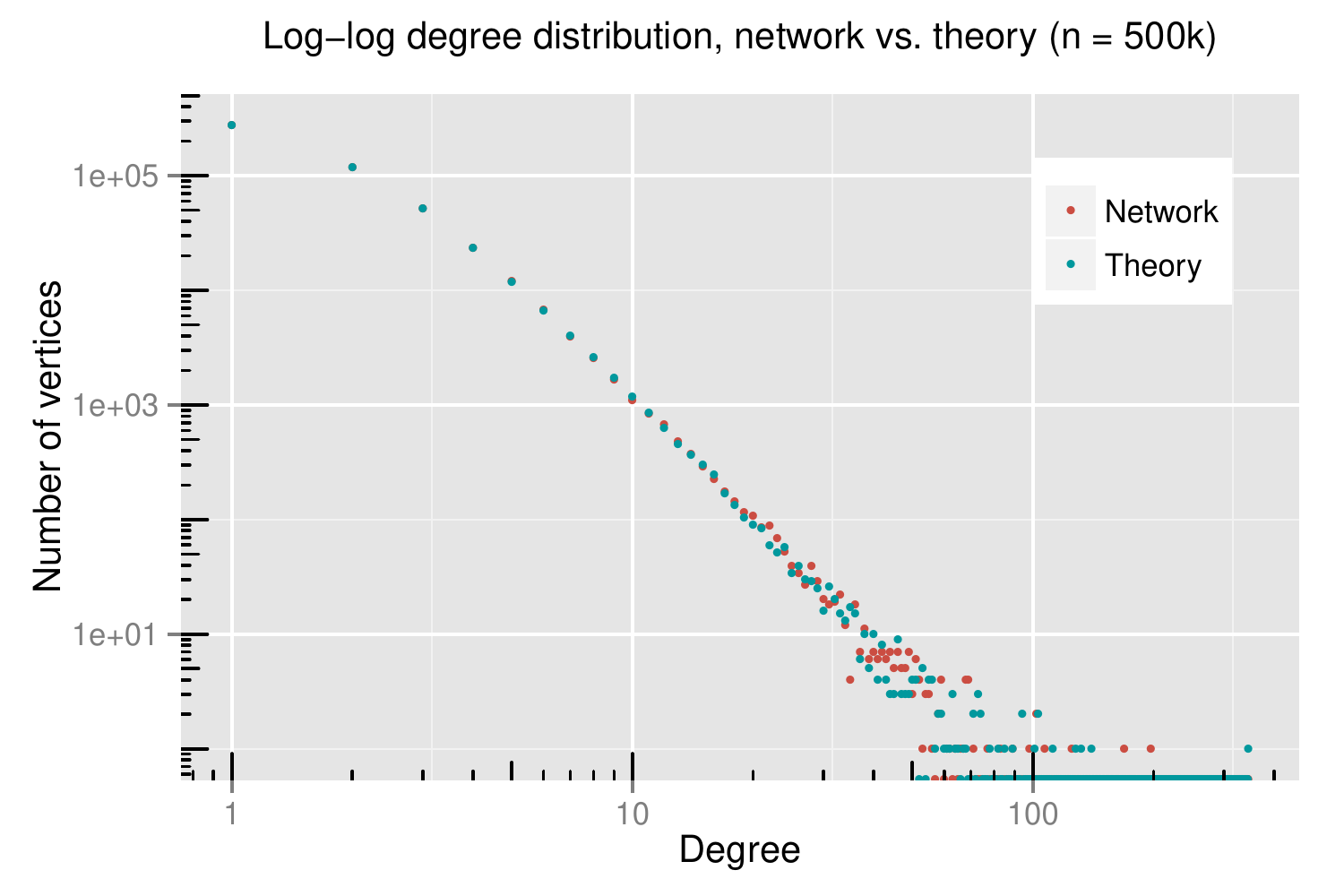}
	\caption{Log log plot showing the limiting degree distribution (red) and simulated network degree distribution (blue) with network size $n=500,000$ and a corresponding sample of the same size from the predicted degree distribution.  The model parameters are taken as $\alpha=6, \beta =1$ and the change point $\gamma=.5$. We discuss other values of the parameters in Section \ref{sec:disc}. }
	\label{fig:log-log}
\end{figure}

\subsection{Change point detection}
\label{sec:res-stat}
The aim of this Section is to formulate a non-parametric estimator for the change point based on observations of the network and state consistency results for this estimator. We first need some notation. 
For fixed $k\geq 1$ let ${N}_n(k,m)$ denote the number of vertices with degree $k$ in the tree $\cT_m$. Rescaling time by $n$, for $0\leq t\leq 1$, let $\hat{N}_n(k,t) = N_n(k,nt)$. Finally define  
\begin{equation}
\label{eqn:deg-time-def}
	\hat{p}^n(k,t) = \frac{\hat{N}_n(k,t)}{nt}, \qquad 0\leq t\leq 1,
\end{equation}
namely the proportion of vertices with degree $k$ at time $nt$. The $k=1$ case corresponds to the number of leaves. To ease notation in the displays below, write $\hat{p}^n(1,t) = \hat{p}^n_t$. Now 
define the continuous function,
	\begin{equation}
	\label{eqn:p-inf-zero-form}
		 p^{\sss(\infty)}_t= \begin{cases}
			\frac{2+\alpha}{3+2\alpha} & \text{if } 0\leq t\leq \gamma\\
			\frac{2+\beta}{3+2\beta}\left(1-\left(\frac{\gamma}{t}\right)^{\frac{3+2\beta}{2+\beta}}\right)+\frac{\gamma}{t}\left(\frac{2+\alpha}{3+2\alpha}\right)\left(\frac{\gamma}{t}\right)^{\frac{1+\beta}{2+\beta}} & \text{if } \gamma \leq t\leq  1. 
		\end{cases}
	\end{equation}
We will prove in Section \ref{sec:leaf-expec} that for each fixed $0< t\leq 1$, $p_t^{\sss(\infty)}$ will represent the limiting proportion of leaves in $\cT_{nt}$. To simplify notation in the sequel, define the function $\delta:\bR_+\to [0,1]$ by the prescription
\begin{equation}
\label{eqn:delta-u-def}
	\delta_u:= \frac{1+u}{2+u}, \qquad u\geq 0.
\end{equation}
Note that $p^{\sss(\infty)}_t = p^{\sss(\infty)}_\gamma$ for $t\leq \gamma$.  Now define the positive function $\set{\sigma_M(t):0\leq t\leq 1}$ via the formulae
\begin{equation}
\label{eqn:simga-t-def}
	\sigma^2_M(t):=\left\{\begin{array}{ll}
		t^{2\delta_\alpha}\left[\delta_{\alpha}p^{\sss(\infty)}_\gamma(1-\delta_{\alpha}p^{\sss(\infty)}_\gamma)\right] & \mbox{ if } 0\leq t\leq \gamma,\\
		& \\
		\gamma^{2\delta_\alpha}\left(\frac{t}{\gamma}\right)^{2\delta_\beta} \delta_\beta p^{\sss(\infty)}_t(1- \delta_\beta p^{\sss(\infty)}_t), & \mbox{if } \gamma < t\leq 1.  
	\end{array}
	\right.
\end{equation}
For later use define the functions 
\begin{equation}
\label{eqn:simga-true-def}
	\sigma^2(t):=\left\{\begin{array}{ll}
		\left[\delta_{\alpha}p^{\sss(\infty)}_\gamma(1-\delta_{\alpha}p^{\sss(\infty)}_\gamma)\right] & \mbox{ if } 0\leq t\leq \gamma,\\
		& \\
		 \delta_\beta p^{\sss(\infty)}_t(1- \delta_\beta p^{\sss(\infty)}_t), & \mbox{if } \gamma < t\leq 1,  
	\end{array}
	\right.
\end{equation}
and 
\begin{equation}
\label{eqn:mut-def}
	\mu(t):= \left\{
	\begin{array}{ll}
		-\frac{\delta_{\alpha}}{t^{\delta_\alpha+1}} & 0< t\leq \gamma\\
		&\\
		-\frac{\delta_{\beta} \gamma^{\delta_{\beta} - \delta_{\alpha}}}{t^{\delta_\beta+1}} & \gamma< t\leq 1\\
	\end{array}
	\right.
\end{equation}

Define the diffusion $\set{M(t):0\leq t\leq 1}$ via the prescription 
\begin{equation}
\label{eqn:dmt-def}
	dM(t) = \sigma_M(t) dB(t), \qquad 0\leq t\leq 1. 
\end{equation}
Here $\set{B(u): u\geq 0}$ is standard Brownian motion on $\bR_+$. 
Thus $M$ is essentially a deterministic time change of $B(\cdot)$ namely 
\begin{equation}
\label{eqn:gt-def}
	\phi(t) = \int_0^t \sigma^2_M(s)ds, \qquad \set{M(t):0\leq t\leq 1} \stackrel{d}{=} \set{B(\phi(t)): 0\leq t\leq 1}.
\end{equation}
In particular $M(\cdot)$ is a Gaussian process on $[0,1]$. 
Finally define the functions 
\begin{equation}
\label{eqn:gtx-def}
	g(t):=\left\{\begin{array}{ll}
		\frac{1}{t^{\delta_{\alpha}}}  & \text{ if } 0< t\leq  \gamma,\\
		& \\
		\frac{\gamma^{\delta_\beta -\delta_\alpha}}{t^{\delta_\beta}} & \text{ if } \gamma < t\leq  1. 
	\end{array} 
	\right.
\end{equation}

Define the process 
\begin{equation}
\label{eqn:gt-process-def}
	G(t)= g(t) M(t), \qquad 0< t\leq 1. 
\end{equation}
By Ito's formula $G(\cdot)$ solves the SDE 
\begin{equation}
\label{eqn:sde-for-G}
	dG(t)= \mu(t) M(t) dt + \sigma(t) dB(t),
\end{equation}
where $\sigma(\cdot)$ and $\mu(\cdot)$  are as in \eqref{eqn:simga-true-def} and \eqref{eqn:mut-def} respectively. Then we have the following result. 

\begin{thm}\label{thm:leaves-fclt}
	Consider the process of re-centered and normalized number of leaves 
	\begin{equation}
	\label{eqn:gnt-def}
		G_n(t):= \frac{\hat{N}_{n}(1,t) - nt p_t^{\sss(\infty)}}{\sqrt{n}}, \qquad 0\leq t\leq 1, 
	\end{equation}
	with linear interpolation between time points. Then $G_n\weakc G$ where $G$ is the diffusion defined in \eqref{eqn:sde-for-G} and convergence is with respect weak convergence on $D([0,1])$ with respect to the usual Skorohod metric. 
\end{thm}

For the rest of this section, let $p_n(m)$ denote the proportion of leaves (degree one vertices) in $\cT_m$. 
Fix $\eps> 0$. We will define two functions on the interval $[\eps,1]$. Let 
\begin{equation}
\label{eqn:bc-def}
	\bcpmf^{\sss(n)} = \frac{1}{n(t-\eps)}\sum_{m=n\eps}^{nt} p_n(m), \qquad \eps\leq t\leq 1.
\end{equation}
Let
\begin{equation}
\label{eqn:ac-def}
	\acpmf^{\sss(n)} = \frac{1}{n(1-t)}\sum_{m=nt+1}^{n} p_n(m), \qquad \eps\leq t\leq 1.
\end{equation}
In words, $\bcpmf^{\sss(n)}$ represents the average proportion of leaves in the process between time $n\eps$ and $nt$ while $\acpmf^{\sss(n)}$ represents the same quantity but after time $nt$. Define the function 
\begin{equation}
\label{eqn:dnt-fn}
	D_n(t):= (1-t)|\bcpmf^{\sss(n)} - \acpmf^{\sss(n)}|, \qquad t\in [\eps,1]. 
\end{equation}
Write $\cM_n$ for the collection of points $t$ for which the corresponding  function value $D_n(t)$ is within $\log{n}/\sqrt{n}$ of the maximum of the function. Precisely, let 
$D_n^* = \max_{t\in [\eps,1]} D_n(t)$ and let 
\begin{equation}
\label{eqn:cmn-def}
	\cM_n:= \set{t\in [\eps,1]: |D_n(t) - D_n^*|\leq \frac{\log{n}}{\sqrt{n}}}.
\end{equation}
Finally let 
\begin{equation}
\label{eqn:estimator-def}
	\hat{\gamma}_n:= \max \set{t: t\in \cM_n}. 
\end{equation}
The functionals $D_n^*, \cM_n,\text{ and } \hat{\gamma}_n$ all depend on $\eps$ but we suppress this dependence to ease exposition below.

\begin{thm}
\label{thm:cp-estimator-cons}
Assume that the change point $\gamma > \eps$. Then the estimator $\hat{\gamma}_n \probc \gamma$ and in fact 
\begin{equation}
\label{eqn:gamma-gammhat}
	|\hat{\gamma}_n - \gamma| = O_P\left(\frac{\log{n}}{\sqrt{n}}\right)
\end{equation}
 Thus $\hat{\gamma}_n$ is a consistent estimator for the change point $\gamma$.
\end{thm}
\begin{rem}
	The $\eps$-truncation away from zero is a technical compensation for the factor $t$ in the denominator in \eqref{eqn:bc-def}. Technically one should be able to choose a sequence $\eps_n\downarrow 0$ slowly enough such that the above result (modified using this sequence $\eps_n$ instead of the fixed $\eps$) is true. This would make the assumption of $\gamma >\eps$ irrelevant in the statement of the Theorem. 
\end{rem}

\begin{rem}
	The threshold $\log{n}/\sqrt{n}$ in \eqref{eqn:cmn-def} was arbitrary in the sense that if we chose the threshold to be $\omega_n/\sqrt{n}$ where $\omega_n\to \infty$ arbitrarily slowly then the corresponding estimator would satisfy \eqref{eqn:gamma-gammhat} with bound $\omega_n/\sqrt{n}$. 
\end{rem}

\begin{figure}[htbp]
	\centering
		\includegraphics[scale=.8]{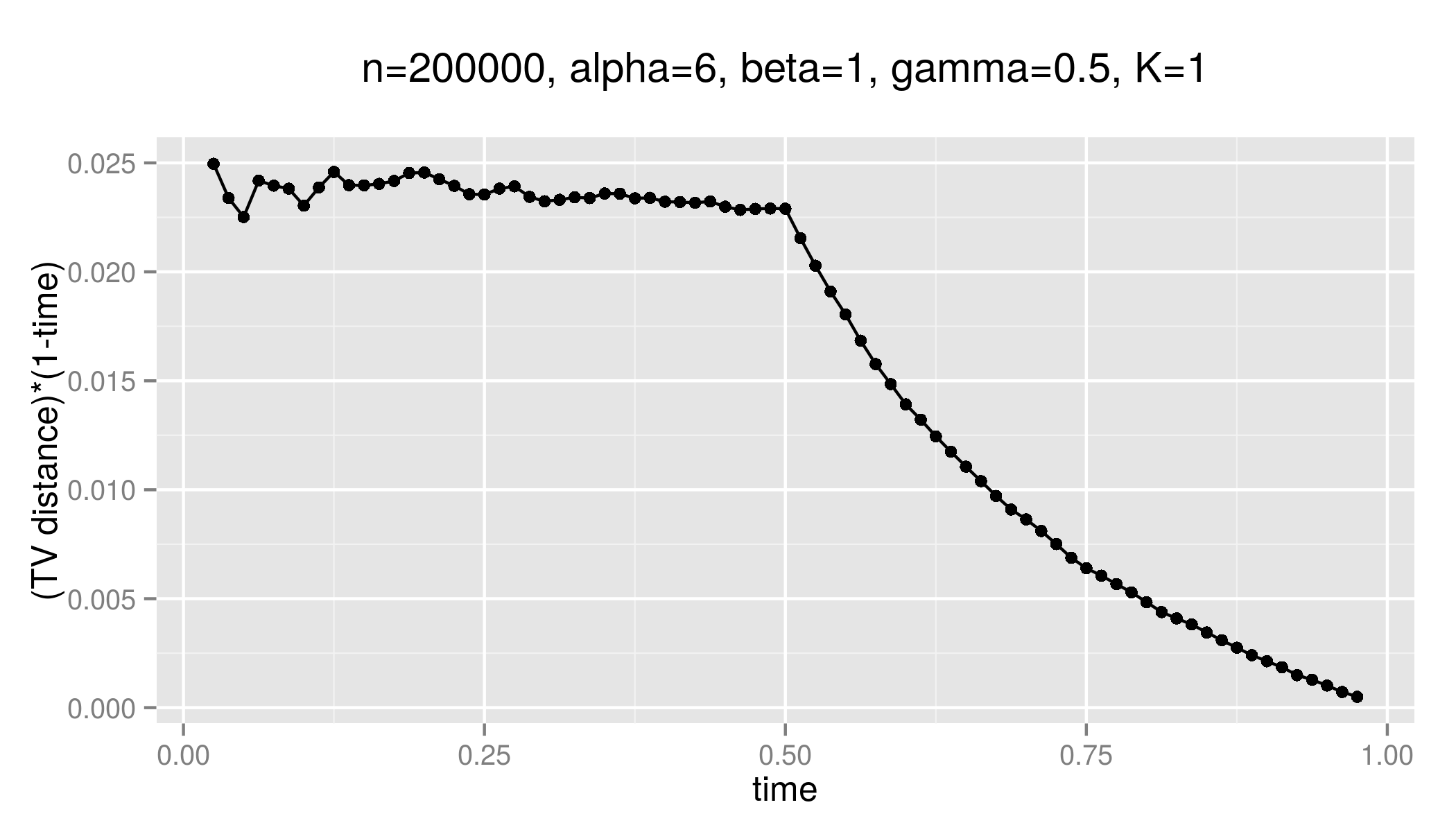}
	\caption{The function $D_n(t)$ with network size $n=200,000$, and model parameters $\alpha=6, \beta =1$ and the change point $\gamma=.5$ as in Figure \ref{fig:log-log}.}
	\label{fig:dtv}
\end{figure}

\begin{rem}
	See Figure \ref{fig:dtv} for a figure based on simulations for the function $D_n(t)$ with $\eps$ taken to be zero. 
\end{rem}

%
%
%

\section{Discussion}
\label{sec:disc}
We now discuss the relevance of our results, their connections to existing literature and possible extensions of the results in this paper.
 
\subsection{\bf Multiple change points}
\label{sec:multiple-cp}
 The proof techniques carry over in a straightforward fashion to the general setting of multiple change points. Fix time points $0< \gamma_1 < \gamma_2 <\cdots \gamma_k < 1$ and parameters $\alpha, (\beta_i)_{1\leq k}$. As before write $\mvtheta = (\alpha, (\beta_i)_{1\leq i\leq k}, (\gamma_i)_{1\leq i \leq k})$ for the parameter set. Consider the random tree $\cT_n = \cT_{\mvtheta,n}$ where 
		\begin{enumeratei}
			\item In the interval $\set{1 < t \leq \gamma_1 n}$, vertices use the attachment scheme driven by $\alpha$ (namely each new vertex attaches to an existing vertex with probability proportional to out-degree  $+1+\alpha$). 
			\item In subsequent intervals $\set{\gamma_{j}n < t \leq \gamma_{j+1} n }$ where $1\leq j\leq k-1$, vertices perform the attachment scheme driven by the parameter $\beta_j$. Here we use the convention $\gamma_0 = 0, \gamma_{k+1} = 1$. 
		\end{enumeratei}  
		As in Section \ref{sec:prelim} define the point processes $\cP_\alpha, \cP_{\beta_i}$ and for fixed $j\geq 1$, the point processes $\cP_\alpha^j, \cP_{\beta_i}^j$. To simplify notation, for any $t\geq 0$ and point process $\cP$, set $\cP[0,t]$ for the \emph{number} of points in the interval $[0,t]$. Define the constants
		\begin{equation}
		\label{eqn:multi-a-pi-ef}
			\pi_j = \gamma_{j+1} - \gamma_{j}, \qquad a_j = \frac{1}{2+\beta_j}\log{\frac{\gamma_{j+1}}{\gamma_j}}. 
		\end{equation}
	Note that $\mvpi=(\pi_0, \pi_1, \ldots, \pi_k)$ is a probability mass function.  Write $\epoch$ for a random variable with distribution $\mvpi$ (i.e. $\pr(\epoch = i) = \pi_i$ for $0\leq i\leq k$). Using the constants $\set{a_i:1\leq i\leq k}$ let  $G_{a_i}$ denote corresponding truncated exponential distributions as in \eqref{eqn:ga-cdf-def} and let $\age_i$ denote a random variable with distribution $G_{a_i}$. Now construct the random variable $\timea$ as follows:
	\begin{enumeratea}
		\item Generate a collection of independent random variables $\epoch$ and $\set{\age_i:1\leq i\leq k}$ with distributions specified as above. 
		\item Conditional on $\epoch = i$, let 
		\[\timea = \age_i + \sum_{j=i+1}^k a_j,\]
		where again by convention, if $\epoch = 0$, $\age_0 =0$ and so $\timea = \sum_{j=1}^k a_i$. 
	\end{enumeratea}
		Construct a positive integer valued random variable $D_{\mvtheta}$ as follows:
		\begin{enumeratei}
			\item Generate $\epoch\sim \mvpi$ as above and the corresponding random variable $\timea$. 
			\item If $\epoch$ takes a non-zero value $1\leq i \leq k$, conditional on $\epoch =i$, generate the switching point process $\cP_{\star}$ on the interval $[0,\timea]$ as follows: 
			\begin{enumeratea}
				\item {\bf Initialization:} In the interval $[0,\age_i]$, start with $\cP_{\star} = \cP_{\beta_i}$.  Suppose by time $\age_i$, $\cP_{\star}[0,\age_i] = k$. Now generate a point process $\cP_{\beta_{i+1}}^{k+1}$ and let $\cP_{\star}[0,\age_i+a_{i+1}] = \cP_{\star}[0,\age_i]+\cP_{\beta_{i+1}}^k[0,a_{i+1}]$.
				\item {\bf Recursion:} For each subsequent interval $[a_j, a_{j+1}]$ with $j> i$, conditional on $\cP_{\star}[0,\age_i+a_{i+1}+\cdots a_j] = k_j $, generate the point process $\cP_{\beta_{j+1}}^{k_j+1}$. Define \[\cP_{\star}[0,\age_i+a_{i+1}+\cdots a_{j+1}] = \cP_{\star}[0,\age_i+a_{i+1}+\cdots a_j] + \cP_{\beta_{j+1}}^{k_j+1}[0,a_{j+1}].  \]
				Iterate until the last interval resulting in $\cP_{\star}[0,\timea]$.
			\end{enumeratea}
			   Now define $D_{\mvtheta} = 1+ \cP_{\star}[0,\timea]$. 
			\item If $\epoch =0$, so that $\timea = a_1+\cdots a_k$, generate a random variable $D_\alpha$ with distribution $\vp_\alpha$ as in \eqref{eqn:pk-alpha-def}. Conditional on $D_\alpha$, generate $\cP_{\star}$ in the interval $[0,a_1]$ with distribution $\cP_{\beta_1}^{D_\alpha}$ and then sequentially proceed as in (ii). In this case, define $D_{\mvtheta} = D_\alpha + \cP_{\star}[0,\timea]$. 
		\end{enumeratei}
		Write $p_{\mvtheta}(\cdot)$ for the pmf of $D_{\mvtheta}$. As before for $k\geq 1$, let $N_n(k)$ denote the number of vertices with degree $k$ in $\cT_n$.
		Then we have the following result. 
		\begin{thm}\label{thm:multi-change-pt}
			As $n\to\infty$ we have 
	 \[\frac{N_n(k)}{n} \probc p_{\mvtheta}(k). \]
	  Further there exist constants $0 < c < c^\prime$ such that for all $k\geq 1$
	 \begin{equation}
	 \label{eqn:mvtheta-multiple-tail}
	 	\frac{c}{k^{\alpha+2}} \leq \pr(D_{\mvtheta} \geq k) \leq  \frac{c^\prime}{k^{\alpha+2}}.
	 \end{equation} 
		\end{thm}
%
		\subsection{\bf Change point detection:} This problem has a vast history owing to its obvious importance in applications in fields ranging from quality control and reliability of industrial processes, in particular quick detection of process failure in production, to fields such as signal processing (e.g. biomedical data including neuronal spike data and seismic data), automatic segmentation of signals into stationary segments via identification of change points etc. While it is impossible to provide a representative sampling of this area,  we direct the interested reader to \cite{basseville1993detection,csorgo1997limit,brodsky1993nonparametric,carlstein1988nonparametric,carlstein1994change,shiryaev1963optimum,shiryaev2007optimal,siegmund1985sequential} and the references therein for an overview of just some of the statistical methodology as well as applications. 
		
		In this context, recall the motivating example of an independent stream of data $\set{X_i:1\leq i\leq n}$ with a change point in the distribution from $F$ to $G$ at time $n\gamma$ described in Section \ref{sec:intro}. Let $_tH^{\sss(n)}(\cdot)$ and $H_t^{\sss(n)}$ denote the empirical distribution of the data before and after $t$ namely 
		\[_tH^{\sss(n)} := \frac{1}{nt} \sum_{i=1}^{nt} \delta_{X_i},\qquad H_t^{\sss(n)}:= \frac{1}{n(1-t)} \sum_{i=nt+1}^{n} \delta_{X_i}, \qquad 0< t< 1. \]
		Now define 
		\[D_n(t):= t(1-t)\dist(_tH^{\sss(n)},H^{\sss(n)}_t),\]
		where $\dist$ is any standard notion of distance between probability distributions on $\bR$ e.g. Kolmogarov-Smirnov supremum norm or total variation distance. Finally define \[\hat{\gamma}_n = \arg\max_{t\in [0,1]} D_n(t).\]
		 Then in \cite{carlstein1988nonparametric} it is shown that $\hat{\gamma}_n$ is a consistent estimator of $\gamma$. This was partial motivation for our estimator. Note the ``asymmetry'' as a function of $t$ between the ``classical'' context and the model with change point highlighting the non-ergodic nature of the evolution of the model after the change point. 
		
		A second point to note is that we use information on leaf densities in the large network $n\to\infty$ limit. As in \cite{resnick2015asymptotic}, one should be able to build on the functional CLT for leaf counts to establish a joint functional CLT for $\set{\hat{N}_n(k,t): 1\leq k \leq K, 0\leq t\leq 1}$ after proper normalization and re-centering for any fixed $K\geq 1$. Modifying the estimator in Section \ref{sec:consistency} should enable one to get estimators that perform better for finite $n$. 

		\subsection{\bf Temporal networks and change points: }  As described in the introduction, the availability of data on real world networks over the last few years has motivated development of mathematical methodology in a wide array of fields including computer science, statistical physics and probability to make sense of this data. With regards to problems philosophically similar to change point detection, analogous to segmentation and boundary detection  \cite{tsybakov1994multidimensional,ma1997edge}, there has been a significant amount of work detecting anomalous subgraphs and motifs {\it within} networks, see e.g. \cite{eberle2007discovering,akoglu2010oddball,noble2003graph}, for a wide-ranging survey see \cite{chandola2009anomaly}. This also includes anomalous edge detection via link prediction algorithms \cite{huang2006link}.   With regards to detection of change points in temporal (time-varying) network data and in particular structural properties of these objects see \cite{sun2007graphscope} that posits an algorithmic approach based on minimum description length to understand evolving communities in social networks. For statistically grounded approaches see \cite{moreno2013network,priebe2005scan,heard2010bayesian,mcculloh2011detecting,zhang2013online,duan2009community,sharpnack2012changepoint}. See \cite{clauset-peel} for an overview of the state of the art regarding change point detection in networks and develops new statistical methodology using a generalized hierarchical random graph model (GHRG) and various likelihood ratio based test statistics to detect existence of change points via online detection algorithms. This paper also studies the performance of these algorithms on simulated as well as real data including the MIT proximity data \cite{eagle2006reality} and the ENRON email network data. See \cite{yudovina2015changepoint,marangoni2015sequential} for rigorous analysis of models where each time slice of the model is assumed to be an Erdos-Renyi random graph. 
\subsection{\bf Preferential attachment:}
	This model has become one of the standard workhorses in the complex networks community, in particular for its ability to give a generative reason for the power law/heavy tailed degree distribution observed in an array of real world systems. At this point it is impossible to compile a representative list of references, we will try to give an overview, restricting ourselves as far as possible to papers close in spirit to this paper; see \cite{szymanski1987nonuniform} where it was introduced in the combinatorics community,  \cite{barabasi1999emergence} for bringing this model to the attention of the networks community,  \cite{newman2003structure},\cite{dorogovtsev2002evolution} for survey level treatments of a wide array of models, \cite{Bollobas:2001:DSS:379831.379835} for the first rigorous results on the asymptotic degree distribution, and \cite{cooper2003general}, \cite{bollobas2003mathematical}, \cite{rudas-2}, and \cite{durrett-rg-book} and the references therein for more general models and results.
	
	We are not aware of other analysis of the effect of change point in structural properties of such network models. There has been a lot of recent interest in understanding and detecting the ``initial seed'' \cite{bubeck-mossel,bubeck-devroye-lugosi,kortchemski}. Here one starts with an initial ``seed graph'' at time $m=0$ and then performs preferential attachment started from that seed. The aim is then to estimate this initial seed based on an observation of the network at some large time $n$. While different from this paper, this body of work again emphasizes the sensitive dependence on initial conditions for such network models.   
	 \subsection{\bf Proof techniques:} \label{sec:disc-proof} A number of techniques have been developed to rigorously analyze functionals such as asymptotic degree distributions (see \cite{durrett-rg-book,van2009random} for nice pedagogical treatment). The standard technique involves writing down recursions for the expected degree distribution $\E(N_n(k))$ using the prescribed dynamics of the process, to show that these expectations (normalized by $n$) converge in the limit and then showing that the deviations $|N_n(k) - \E(N_n(k))|$ are small via concentration inequalities.
	 
	 In this paper, for understanding structural properties we use a different technique, essentially embedding the discrete time model in a corresponding ``continuous time'' branching process $\set{\BP_{\mvtheta}^n(t): t\geq 0}$ (based on the Athreya-Karlin embedding of urn processes \cite{athreya1968}). This explains the various point processes that arise in the description of the limiting degree distribution. While mathematically more involved, this technique gives more insight into the results as it elucidates the natural time scale of the process. In various other settings this technique has resulted in the study of much more general functionals of the process such as the spectral distribution of the adjacency matrix \cite{bhamidi2012spectra} and has been used to derive asymptotic results in ``non-local'' preferential attachment models \cite{bhamidi2012twitter}. In this paper the technique also allows one to intuitively understand why the degree exponent does not change. We advise the reader to come back to the text below after going through the proofs but let us explain the basic intuition here.  In the continuous time version, the process grows exponentially and in particular takes time $\tau_{\gamma n}\approx \frac{1}{2+\alpha}\log{\gamma n}+O_P(1)$ to get to size $n\gamma$. At this time there is a change in the evolution where each vertex adopts attachment dynamics driven by the parameter $\beta$. However owing to the exponential growth rate, the time for the process to get to size $n$ is $\tau_n \approx \tau_{\gamma n} + a$ where $a$ is as in \eqref{eqn:a-def}. Thus the process does not have enough time for the dynamics with attachment parameter $\beta$ to change the degree exponent (since we only have to wait an $O(1)$ extra units of time to get to system size $n$ from $\gamma n$). These ideas are made mathematically rigorous in the next few sections. For the interested reader, much of the foundational work on continuous time branching processes relevant for this paper can be found in \cite{jagers-nerman-1,jagers-nerman-2, jagers1975branching}.

	\subsection{Empirical dependence of the convergence on parameter values:}
	\label{sec:subsec-empirical}
	
Recall that the Gaussian process defined in \eqref{eqn:sde-for-G} underlying the main consistency result Theorem \ref{thm:cp-estimator-cons} depends on $\mvtheta=(\alpha,\beta,\gamma)$. One consequence of this dependence is that when the parameter values $\alpha$ and $\beta$ are close, the change point becomes harder to detect in the sense that larger $n$ is required to get good estimates. This is most easily seen in terms of the fluctations of the proportion of leaves in the graph. 

\begin{figure}[htbp]
	\centering
		\includegraphics[scale=0.75]{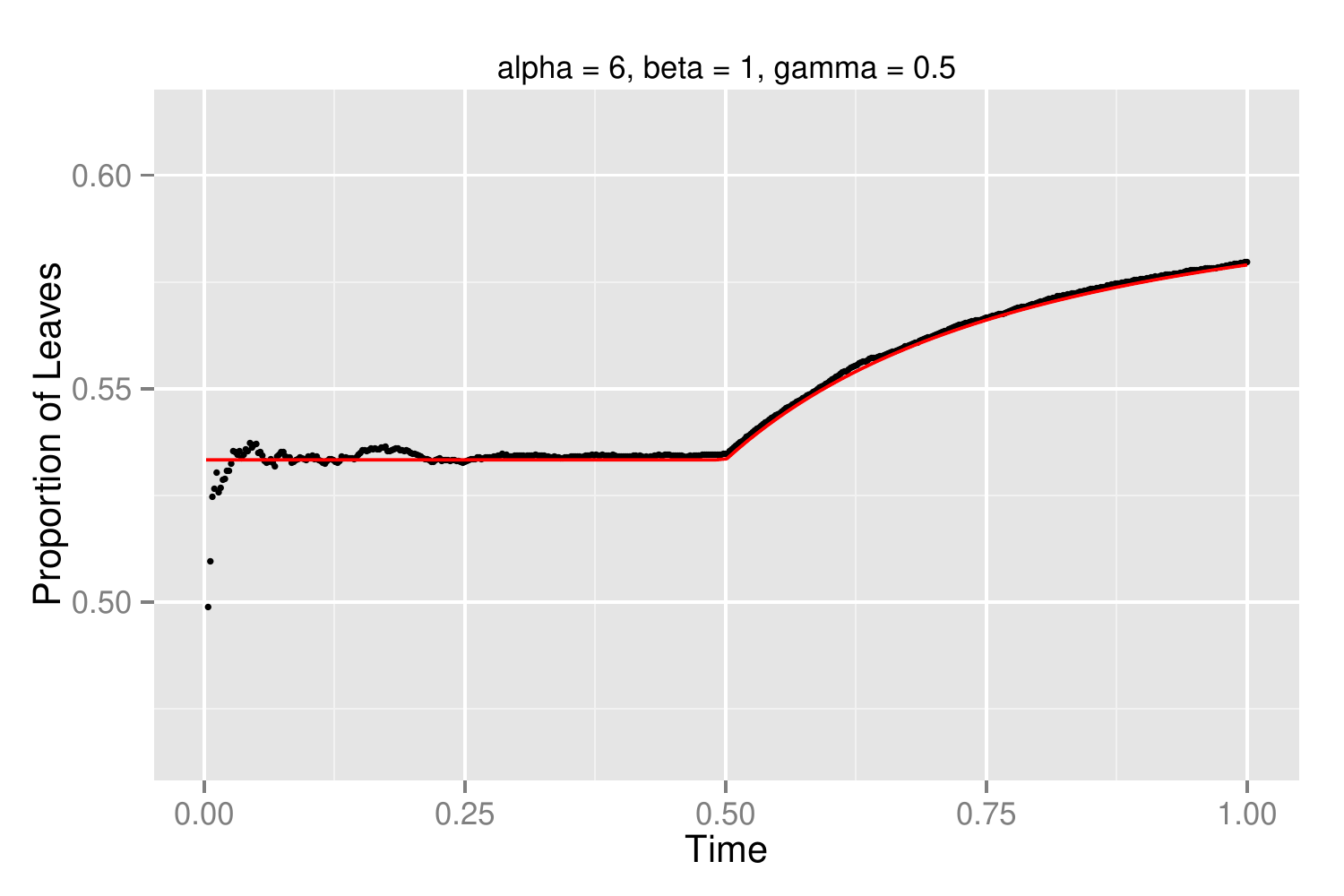}
	\caption{Empirical proportion of leaves in a simulation with $n=200,000, \alpha=6, \beta=1, \gamma=0.5$. The red line represents the theoretical predictions in \eqref{eqn:p-inf-zero-form}. }
	\label{fig:leaves-good}
\end{figure}

\begin{figure}[htbp]
	\centering
		\includegraphics[scale=0.75]{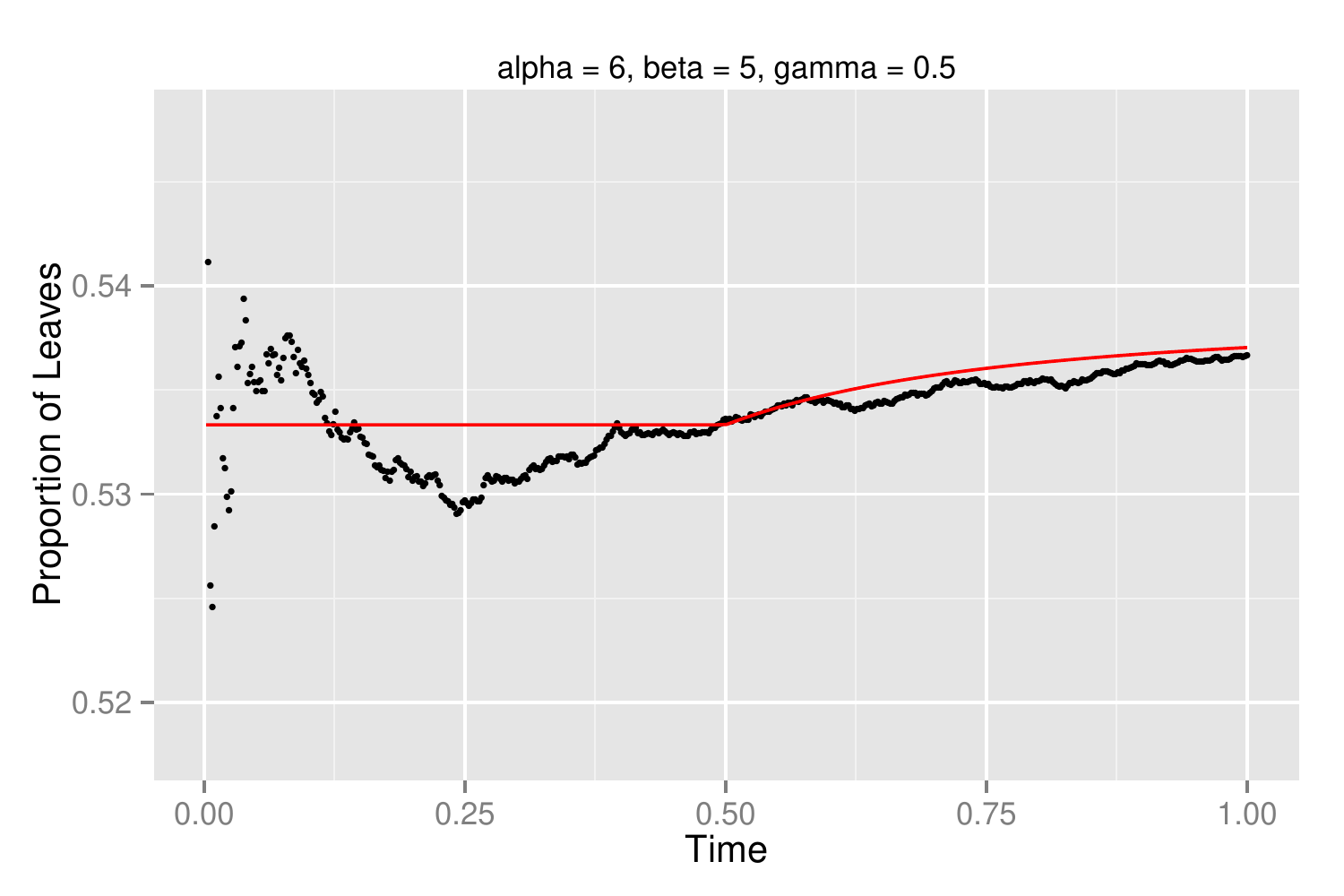}
	\caption{Empirical proportion of leaves in a simulation with $n=200,000, \alpha=6, \beta=5, \gamma=0.5$. The red line represents the theoretical predictions in \eqref{eqn:p-inf-zero-form}. }
	\label{fig:leaves-bad}
\end{figure}
	
In both Figures \ref{fig:leaves-good} and \ref{fig:leaves-bad}, the preferential attachment process starts with $\alpha=6$ and decreases, to $\beta=1$ in \ref{fig:leaves-good} and $\beta=5$ in \ref{fig:leaves-bad}. Furthermore the predicted behavior (red line) is almost the same: the proportion of leaves is constant up to the change point $\gamma=0.5$ and then increases, consistent with a \textit{decrease} in the attachment parameter.

Despite the sizes of the final graphs in both simulations being $n=200,000$ vertices, at first glance the fluctuations appear much greater in the latter case. On closer examination however, this is simply an illusion of the axes. In essence, when the shift in parameters is smaller, the change in the proportion of leaves pre- and post-$\gamma$ is smaller compared to the natural fluctuations in the proportion of leaves which is of order $\sqrt{n}$ (Theorem \ref{thm:leaves-fclt}). Therefore any difference is more difficult to detect for same $n$. This is not surprising, but worth noting in practice.

\section{Proofs}
\label{sec:proofs}

As described in Section \ref{sec:disc-proof}, the main conceptual idea is a continuous time embedding of the discrete time process. 
We start in Section \ref{sec:proof-idea} by describing this embedding and deriving simple properties. Then in Section \ref{sec:convg-degree} we prove Theorem \ref{thm:deg-dist}. Section \ref{sec:proof-tail-bd} proves the assertion that the degree exponent does not change. Section \ref{sec:max-degree-proof} analyzes asymptotics for the maximal degrees. Section \ref{sec:leafclt-proofs} contains an in-depth analysis of the density of leaves and proves Theorem \ref{thm:leaves-fclt}. Section \ref{sec:consistency} then uses this Theorem to prove the consistency of the estimator namely Theorem \ref{thm:cp-estimator-cons}.

\subsection{Preliminaries}
\label{sec:proof-idea}

We start with the following definition. To ease notation, for the rest of the paper we use $\gamma n$ instead of $\lfloor \gamma n \rfloor$. 

\begin{defn}[Continuous time branching process]
	Fix $\alpha > 0$. We let $\set{\BP_\alpha(t): t\geq 0}$ be a continuous time branching process driven by the point process $\cP_\alpha$ defined in \eqref{eqn:cp-alpha-def}. Precisely:
	\begin{enumeratea}
		\item At time $t=0$ we start with one individual called the root $\rho$ with an offspring point process with distribution $\cP_\alpha^\rho \stackrel{d}{=} \cP_\alpha$. The times of this point process represent times of birth of new offspring of $\rho$. 
		\item Every new vertex $v$ that is born into the system is given its own offspring point process $\cP_\alpha^v \stackrel{d}{=} \cP_\alpha$, independent across vertices. 
	\end{enumeratea}
\end{defn}

\begin{figure}
\centering
\begin{minipage}{.5\textwidth}
  \centering
  \includegraphics[height=55mm, width=80mm]{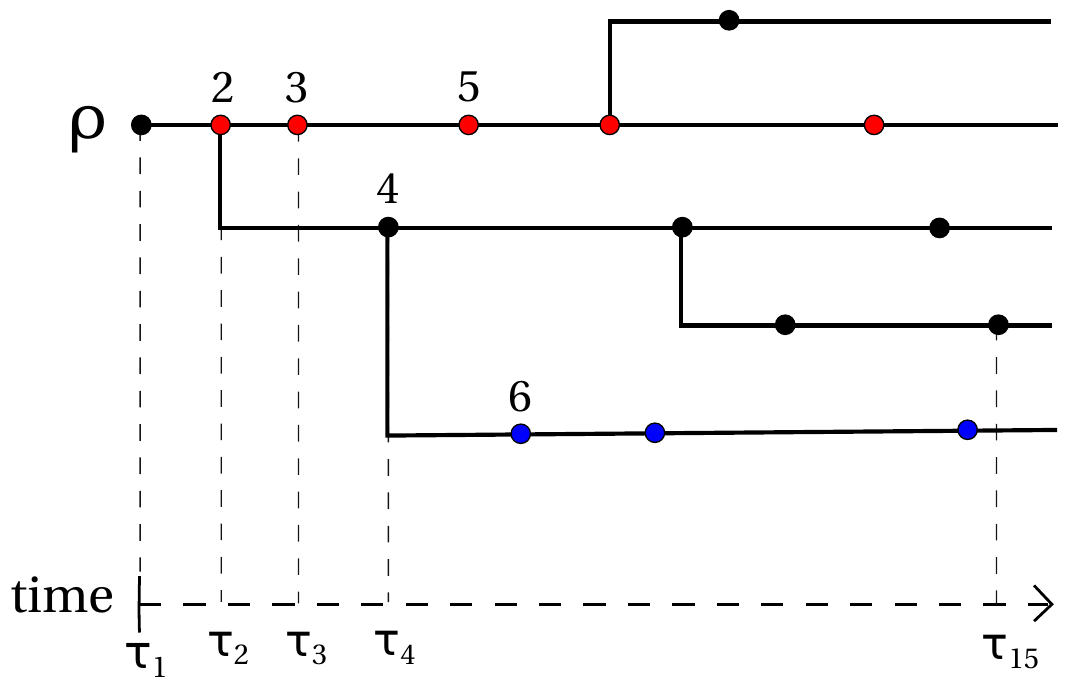}
  \captionof{figure}{The process $\BP_{\alpha}(\cdot)$ in continuous time starting from the root $\rho$ and stopped at $\tau_{15}$.}
  \label{fig:test1}
\end{minipage}%
\begin{minipage}{.5\textwidth}
  \centering
  \includegraphics[height=50mm, width=70mm]{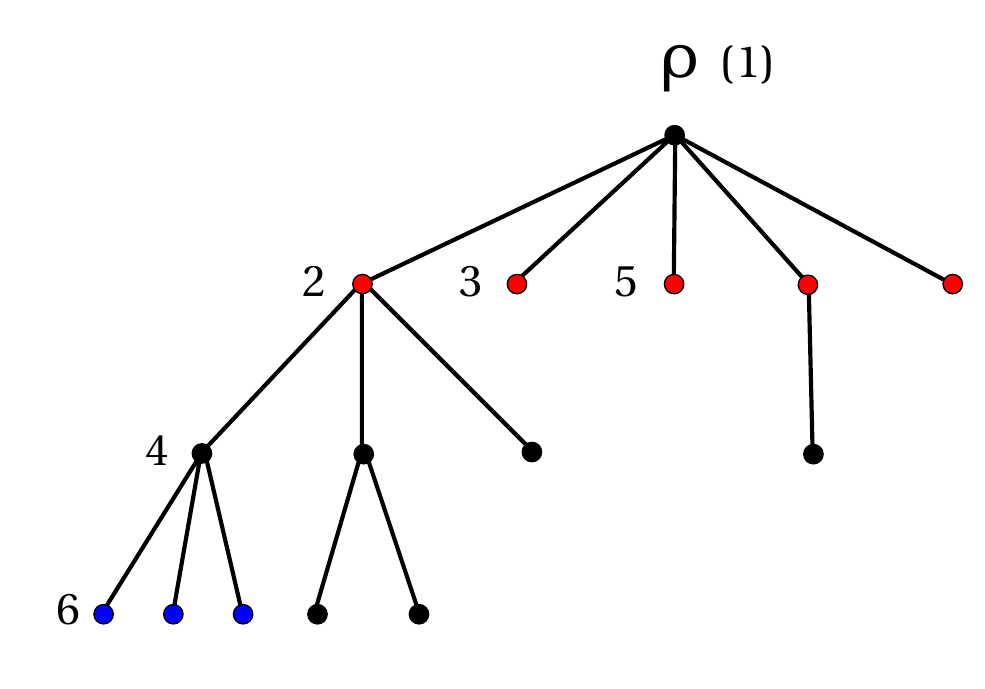}
  \captionof{figure}{The corresponding discrete tree containing only the genealogical information of vertices in $\BP_{\alpha}(\tau_{15})$.}
  \label{fig:test2}
\end{minipage}
\end{figure}

Label vertices using integer labels according to the order in which they enter $\BP_{\alpha}$ so that the root is labelled as $1$, the next vertex to be born labeled by $2$ and so on.  For fixed $t\geq 0$, we will view $\BP_\alpha(t)$ as a (random) labelled tree representing the genealogical relationships between all individuals in the population present at time $t$. See Figures \ref{fig:test1} and \ref{fig:test2}.  Write $|\BP_{\alpha}(t)|$ for the number of individuals in the tree by time $t$. 
Fix $m\geq 1$ and define the stopping time 
\begin{equation}
\label{eqn:stop-time-def}
	\tau_m:= \inf\set{t: |\BP_{\alpha}(t)| =m}.
\end{equation}
 Since there are no deaths and each individual reproduces at rate at least $1+\alpha$, the stopping times $\tau_m < \infty$ a.s. for all $m\geq 1$.  Now consider the original preferential attachment model where there is no change point. Using properties of the exponential distribution, the following Lemma is easy to check and is just a special case of the famous Athreya-Karlin embedding \cite{athreya1968}. 
\begin{lem}
\label{lem:embed-distr}
Viewed as random rooted trees on vertex set $[n]$ one has  $\BP_\alpha(\tau_n)\stackrel{d}{=} \cT_n$.  In fact the two processes of growing random trees have the same distribution namely 
\[\set{\BP_\alpha(\tau_n): n\geq 1} \stackrel{d}{=} \set{\cT_n: n\geq 1}. \]
\end{lem}

To construct the variant $\cT_n$ where one has a change point, we run $\BP_\alpha(\cdot)$ till time $\tau_{\gamma n}$ (when the original process reaches size $\gamma n$) and then every vertex changes the way it reproduces. More precisely, after this stopping time, an individual with $k$ children would have reproduced at rate $k+1+\alpha$ in the original model but in the change point model this vertex reproduces at rate $k+1+\beta$ and uses the parameter $\beta$ instead of $\alpha$ for each subsequent offspring times. Each new vertex $v$ produced after time $\tau_{\gamma n}$ reproduces according to an independent copy of the point process $\cP_\beta$.   Call the resulting process $\BP_{\mvtheta}^n(\cdot)$ and run the process till time $\tau_n$ when the continuous time process has $n$ individuals. Analogous to \eqref{eqn:stop-time-def}, define the collection of stopping times $\set{\tau_m: 1\leq m \leq n}$ by replacing $\BP_\alpha$ with $\BP_{\mvtheta}^n$. The following is a simple extension of the previous Lemma. 
\begin{lem}\label{lem:embed-distr-cp}
Recall the family of random trees $\set{\cT_{\mvtheta,m}: 1\leq m\leq n}$ generated using the change point preferential attachment model in Section \ref{sec:model-cp-def}. Then, 
\[\set{\BP_{\mvtheta}^n(\tau_m): 1\leq m\leq n} \stackrel{d}{=} \set{\cT_{\mvtheta,m}: 1\leq m \leq n}. \]
\end{lem} 

\begin{rem}
	Note that the processes $\set{\cT_{\mvtheta,m}: 1\leq m \leq n}$ when one has a change point are {\bf not} nested in a nice manner as growing trees for {\bf different values} of $n$. Compare this with the original model (without change point) where we can view the entire sequence $\set{\cT_{n}:n\geq 1}$ as an increasing family of random trees. In the above construction it will be convenient to couple the processes across different $n$ by using a {\bf single} common branching process $\BP_\alpha$ to generate the tree before the change point $\tau_{\gamma n}$ and then let the process evolve independently after the change point for different $n$ using the prescribed dynamics modulated by the attachment parameter $\beta$. Further it will be convenient to allow the process $\BP_{\mvtheta}^n$ to continue to grow after time $\tau_n$ as opposed to stopping it exactly at time $\tau_n$.  
\end{rem}

For future reference, for each vertex $v$, we will use $T_v$ for the time of birth of this vertex into the system. For fixed time $t$ and a vertex $v$ born before time $t$ (namely $T_v\leq t$), we write $d_v(t)$ for the number of children of this vertex by time $t$. Note that for all $v\neq \rho\in \BP_{\mvtheta}^n(t)$, the full degree of $v$ by time $t$ is $d_v(t)+1$. 

 We will need some simple stochastic calculus calculations below to derive Martingales related to processes of interest. Given a process $\set{Z(t):t\geq 0}$ adapted to a filtration $\set{\cF(t):t\geq 0}$, we write $\E(dZ(t)|\cF(t)) = a(t) dt$ for an adapted process $a(\cdot)$ if $Z(t)-\int_0^t a(s) ds$ is a (local) martingale. Similarly write $\var(dZ(t)|\cF(t)) = b(t) dt$ if the process
\[V(t):= \left(Z(t) - \int_0^t a(s)ds\right)^2 - \int_0^t b(s) ds, \qquad t\geq 0,\]
is a local martingale. 

 Now recall that $\BP_\alpha(\tau_{\gamma n})$ is the random tree before the change point. These random trees are distributed as the original preferential attachment model without change point using attachment dynamics with parameter $\alpha$. Using \eqref{eqn:pk-alpha-def} and recalling that $N_n(k,\gamma n)$  denotes the number of vertices with degree $k$ results in the following. 

\begin{lemma}\label{lem:as-pre-convg}
	For each fixed $k\geq 1$ we have 
	${N_{n}}(k,\gamma n)/{\gamma n} \convas p_{\alpha}(k)$, as $n\to\infty$
	where $p_\alpha(\cdot)$ is the probability mass function in \eqref{eqn:pk-alpha-def}. 
\end{lemma}
 Recall that the branching process $\BP_\alpha$ is driven by the offspring point process $\cP_{\alpha}$ and $\cP_{\alpha}(t):=\cP_{\alpha}[0,t]$ is the number of points in $[0,t]$.   Define the process 
\begin{equation}
\label{eqn:malpha-def}
	M_{\alpha}(t):= e^{-t} \cP_{\alpha}(t) - (1+\alpha)(1-e^{-t}), \qquad t\geq 0
\end{equation}
\begin{lemma}\label{lem:palpha-martingale}
	The process $\set{M_{\alpha}(t): t\geq 0}$ is a martingale with respect to the natural filtration of $\cP_{\alpha}$. In particular 
	\begin{equation}
	\label{eqn:mean-cpa}
		\E(\cP_{\alpha}(t)) = (1+\alpha)(e^{t}-1)
	\end{equation}
\end{lemma}
\noindent {\bf Proof:} Write $\set{\cF(t):t\geq 0}$ for the natural filtration of the process. It is enough to show for all $t\geq 0$, $\E(dM_{\alpha}(t)|\cF(t)) = 0$. By construction 
\[\E(d\cP_\alpha(t)|\cF(t))= (1+\alpha+\cP_{\alpha}(t))dt.\] 
Further 
\[\E(dM_\alpha(t)|\cF(t))= e^{-t}\E(d\cP_{\alpha}(t)|\cF(t)) -e^{-t}\cP_{\alpha}(t)dt+(1+\alpha)e^{-t}dt.\]
Elementary algebra completes the proof. The final assertion regarding \eqref{eqn:mean-cpa} follows using the Martingale property of $M_\alpha$ and the initial condition $\cP_{\alpha}(0)=0$. 
\qed

The starting point in the analysis of continuous time branching processes is the so called Malthusian rate of growth parameter $\lambda > 0$ which solves the equation 
\begin{equation}
\label{eqn:malthus-solve}
	\int_0^\infty \lambda e^{-\lambda t}\E(\cP_{\alpha}(t)) dt =1 
\end{equation}
Using Lemma \ref{lem:palpha-martingale} now implies 
\begin{equation}
\label{eqn:malthus-param}
	\lambda = 2+\alpha.
\end{equation}
Let $T_\lambda$ be an exponential random variable with parameter $\lambda$ independent of $\cP_{\alpha}$ and consider the integer valued random variable $\cP_\alpha(T_\lambda)$. Note that \eqref{eqn:malthus-solve} is equivalent to $\E(\cP_\alpha(T_\lambda)) = 1$.  Recall that $D_\alpha$ is a random variable with the (non-change point) degree distribution \eqref{eqn:pk-alpha-def}. It is easy to check that 
$D_\alpha -1 \stackrel{d}{=} \cP_{\alpha}(T_\lambda)$. 
 In particular for $\alpha \geq 0$,
\[\E(\cP_{\alpha}(T_\lambda)\log^+ \cP_{\alpha}) < \infty.\]
Using standard Jagers-Nerman stable age-distribution theory for branching processes \cite{jagers-nerman-1,jagers-nerman-2} now implies the following. 

\begin{prop}
	\label{prop:size-bp-alpha}
	There exists an integrable a.s. positive random variable $W_\alpha$ such that 
	\[e^{-(2+\alpha) t}|\BP_{\alpha}(t)| \stackrel{a.e., \bL^1}{\longrightarrow} W_\alpha.\]
	In particular 
	\begin{equation}
	\label{eqn:tgamma-asymp}
		\tau_{\gamma n} - \frac{1}{2+\alpha} \log{n} \convas W^\prime_{\alpha},
	\end{equation}
	for a finite random variable $W_\alpha^\prime$. 
	
\end{prop}

We conclude this Section with asymptotics for the amount of `` continuous time'' where the attachment dynamics using $\beta$ is valid, namely $\tau_n-\tau_{\gamma n}$. Recall the constant $a$ from \eqref{eqn:a-def}. We will also write $\set{\cF_n(t): t\geq 0}$ for the natural filtration of the process $\set{\BP_{\mvtheta}^n(t): t\geq 0}$. 

\begin{lemma}\label{lem:a-def-asymp}
	Let $\Upsilon_n =\tau_n -\tau_{\gamma n} $ denote the time after the change point in the continuous time embedding. Then 
	\[\sqrt{n}(\Upsilon_n - a) \weakc \frac{1}{2+\beta}\sqrt{\frac{1-\gamma}{\gamma}} Z, \]
	as $n\to \infty$. Here $Z$ is a standard normal random variable. 
\end{lemma}
\noindent{\bf Proof:} Note that $\BP_{\mvtheta}^n(\cdot)$ is a Markov process. Further for $t\geq \tau_{\gamma n}$ conditional on $\BP_{\mvtheta}^n(t)$, the rate at which a new individual is born into the system is given by 
\begin{align}
\lambda(t)&:=\sum_{v\in \BP_{\mvtheta}^n(t)} (d_v(t)+1+\beta)\notag\\
&= (2+\beta)|\BP_{\mvtheta}^n(t)| -1, \label{eqn:rate-descp}
\end{align} 
In particular
\begin{equation}
\label{eqn:upsilon-n-def}
	\Upsilon_n \stackrel{d}{=} \sum_{j=\gamma n}^{n-1} \frac{E_i}{(2+\beta)j-1},
\end{equation}
where $\set{E_i: i\geq 1}$ is a sequence of \emph{iid} rate one exponential random variables. Using Lyapunov's central limit theorem now completes the proof. 
\qed 

Using the distributional characterization in \eqref{eqn:upsilon-n-def} and standard concentration inequalities for sums of independent random variables, one can show the following tail bound on $\Upsilon_n$. We omit the proof.

\begin{lemma}\label{lem:upsion-tail}
	For any $\kappa > 0$ there exists $N = N(\kappa)<\infty$ such that for all $n> N(\kappa)$, 
	\[\pr\left(|\Upsilon_n -a|> \frac{1}{n^{1/3}} \right)\leq \frac{1}{n^{\kappa}}.\]
	In particular by Borel-Cantelli, 
	$\pr\left(|\Upsilon_n -a|\leq  {n^{-1/3}} \text{ eventually }\right)=1. $
\end{lemma}
Here the bound $n^{-1/3}$ was arbitrary. An upper bound of $n^{-(1/2-\delta)}$ with any $\delta >0$ would result in identical result as above. We fix $n^{-1/3}$ for definiteness.  We end this Section by defining the Yule process. Properties of this process will be needed in the next few Sections.

\begin{defn}[Rate $\nu$ Yule process]
\label{def:yule-process}
	Fix $\nu > 0$. A rate $\nu$ Yule process is a pure birth process $\set{Y_\nu(t):t\geq 0}$ with $Y_\nu(0)=1$ and where the rate of birth of new individuals is proportional to size of the current population. More precisely 
	\[\pr(Y_\nu(t+) - Y_\nu(t)|\cF(t)):= \nu Y_\nu(t) dt + o(dt),\]
	where $\set{\cF(t):t\geq 0}$ is the natural filtration of the process. 
\end{defn}

The following is a standard property of the Yule process, see e.g. \cite[Section 2.5]{norris-mc-book}.
\begin{lem}
\label{lem:yule-prop}
Fix time $t >0$ and rate $\nu > 0$. Then the random variable $Y_\nu(t)$, namely the number of individuals in the population by time $t$ has a Geometric distribution with parameter $p=e^{-\nu t}$ namely 
\[\pr(Y_\nu(t)  = k) = e^{-\nu t}(1-e^{-\nu t})^{k-1}, \qquad k\geq 1.\] 	
\end{lem} 

\subsection{Convergence of the degree distribution} 
\label{sec:convg-degree}
In this Section we will prove Theorem \ref{thm:deg-dist}. Recall the description of the limit random variable $D_{\mvtheta}$ in Section \ref{sec:prelim}. It will be easier to deal with the random variable $D_{\mvtheta}^{\out}:= D_{\mvtheta}-1$. Then the distribution of $D_{\mvtheta}^{\out}$ can be written succinctly as: 
\begin{enumeratea}
	\item with probability $\gamma$, $D_{\mvtheta}^{\out}:= Y_{\bc}$ where $Y_{\bc}= D_{\alpha}-1+ N_{\beta}^{D_{\alpha}}[0,a]$;
	\item with probability $1-\gamma$,  $D_{\mvtheta}^{\out} = Y_{\ac}$ where $Y_{\ac}:=X_{\ac}$ and $X_{\ac}$ is as defined in Section \ref{sec:prelim}. 
\end{enumeratea}
Now recall that for any time $t$ and vertex $v$ born before time $t$, $d_v(t)$ denotes the number of children (out-degree) of vertex $v$ at time $t$. For fixed $k\geq 0$ define 
\begin{equation}
\label{eqn:nn-pre-def}
	\bar{N}_n^{\bc}(k):= \sum_{v\in \BP_{\mvtheta}(\tau_n)}\ind\set{T_v \leq \tau_{\gamma n}, d_v(\tau_{n}) \geq k},
\end{equation}
and 
\begin{equation}
\label{eqn:nn-post-def}
	\bar{N}_n^{\ac}(k):= \sum_{v\in \BP_{\mvtheta}(\tau_n)}\ind\set{T_v > \tau_{\gamma n}, d_v(\tau_{n}) \geq k}.
\end{equation}
In words, $\bar{N}_n^{\bc}(k)$ are the number of vertices that were born before the change point and have out-degree at least $k$ by time $\tau_n$ (thus in the tree $\cT_{\mvtheta,n}$) whilst $\bar{N}_n^{\ac}(k)$ is defined analogously but for vertices born after the change point $\tau_{\gamma n}$. The following proposition is equivalent to Theorem \ref{thm:deg-dist}. 

\begin{prop}
	\label{prop:pre-post-convg}
	Fix $k\geq 0$. Then we have 
	\begin{equation}
	\label{eqn:nn-pre-post-convg}
		\frac{\bar{N}_n^{\bc}(k)}{n} \probc \gamma \pr(Y_{\bc}\geq k), \qquad \frac{\bar{N}_n^{\ac}(k)}{n} \probc (1-\gamma) \pr(Y_{\ac}\geq k),
	\end{equation}
	as $n\to \infty$. 
\end{prop} 
The rest of this Section deals with proving this Proposition. 
\subsubsection{\bf Analysis of $\bar{N}_n^{\bc}(\cdot)$ :}
\label{sec:proof-nnbc}
We start with the easier case. We will need some more notation. For fixed $0\leq j, k$, define $\bar{N}_n^{\bc}(j:k)$ for the number of vertices that were born before the the change point $\tau_{\gamma n}$ with out-degree exactly $j$ at time $\tau_{\gamma n}$ that end up with at least $k$ children by time $\tau_n$. Note that
\[\sum_{j\geq k} \bar{N}_n^{\bc}(j:k) = N_{n}(k+1,\gamma n) \]
 namely the number of vertices with total degree $k+1$ (thus out-degree $k$) in the tree before change point $\cT_{\gamma n}$. Recall that Lemma \ref{lem:as-pre-convg}, the asymptotic degree distribution of $\cT_{\gamma n}$ is $D_{\alpha}$ and thus the asymptotic out-degree distribution of the tree $\cT_{\gamma n}$ is $D_{\alpha}^{\out} = D_{\alpha}-1$. Using the form of $Y_{\bc}$, it is thus enough to show for each fixed $0\leq j\leq k$, 
\begin{equation}
\label{eqn:nn-jk-convg}
	\frac{\bar{N}_n^{\bc}(j:k)}{n} \convas \gamma \pr(D_{\alpha}^{\out} = j) \pr(\cP_\beta^{j+1}[0,a]\geq k-j). 
\end{equation} 
We start with the following simple Lemma. 
\begin{lemma}\label{lem:binom-convg}
	Fix $0< p,q< 1$, a sequence of non-negative integer valued random variables $\set{N_n:n\geq 1}$ and a sequence $\set{q_n: n\geq 1} \in [0,1]$. Conditional on $N_n$, let $S_n$ be a Binomial$(N_n,q_n)$ random variable. Further suppose 
	\[\frac{N_n}{n} \convas p, \qquad q_n \to q.\]
	Then $S_n/n\convas pq$.   
\end{lemma}
\noindent{\bf Proof:} We assume we work on a rich enough probability space where we can couple $\set{S_n:n\geq 1}$ with a sequence $\set{\tilde{S}_n:n\geq 1}$ where $\tilde{S_n}$ is Binomial$(np, q_n)$ such that $|S_n - \tilde{S_n}|\leq |N_n-np|$. Standard exponential tail bounds for the Binomial distribution coupled with Borel Cantelli and the hypothesis of the Lemma imply that $\tilde{S}_n/n\convas pq$. Since $|S_n - \tilde{S_n}|/n\leq |N_n/n-p|$, again using the hypothesis of the Lemma completes the proof. \qed

We now proceed with the proof. Analogous to $\bar{N}_n^{\bc}(j:k)$, for each $s\geq 0$ define $\bar{Z}_n^{\bc}((j:k),s)$ for the number of vertices born before the change point $\tau_{\gamma n}$ such that at $\tau_{\gamma n}$ they have out-degree exactly $j$ and further by time $\tau_{\gamma n} +s$ they have degree at least $k$. Then note that conditional on the information at time $\tau_{\gamma n}$, 
\begin{equation}
\label{eqn:zn-binomial}
	\bar{Z}_n^{\bc}((j:k),s) \stackrel{d}{=} \text{Bin}(N_n(j+1, \gamma n),\pr(\cP_\beta^{j+1}[0,s]\geq k-j) )
\end{equation}
Further the random variables of interest $\bar{N}_n^{\bc}(j:k) = \bar{Z}_n^{\bc}((j:k),\Upsilon_n) $ where $\Upsilon_n$ is as in Lemma \ref{lem:a-def-asymp}. Thus writing $a_n^+ = a+n^{-1/3} $ and $a_n^{-} = a-n^{-1/3}$ and using Lemma \ref{lem:upsion-tail},
\begin{equation}
\label{eqn:nn-zn-bound}
	\bar{Z}_n^{\bc}((j:k),a_n^-) \leq  \bar{N}_n^{\bc}(j:k) \leq \bar{Z}_n^{\bc}((j:k),a_n^+) \text{ eventually a.s.}
\end{equation}
Using the Binomial convergence Lemma \ref{lem:binom-convg} and noting that by Lemma \ref{lem:as-pre-convg} and choice of $a_n^+, a_n^-$, the hypothesis of this Lemma are satisfied, implies that 
\[\frac{\bar{Z}_n^{\bc}((j:k),a_n) }{n} \convas \gamma \pr(D_{\alpha}^{\out} = j) \pr(\cP_\beta^{j+1}[0,a]\geq k-j),   \]
where take $a_n$ as either $a_n^+$ or $a_n^-$. Now using \eqref{eqn:nn-zn-bound} proves \eqref{eqn:nn-jk-convg}. This completes the analysis of $\bar{N}_n^{\bc}(\cdot)$. 

\qed

\subsubsection{\bf Analysis of $\bar{N}_n^{\ac}(\cdot)$ :}
We start by setting up some notation. Fix $k\geq 0$ and define the function 
\begin{equation}
\label{eqn:gku-def}
	g_k(u):= \pr(\cP_{\beta}[0,u]\geq k), \qquad u\geq 0. 
\end{equation}
Here $\cP_\beta$ is the offspring point process with attachment parameter $\beta$. Then writing out the form of the distribution of $Y_{\ac}$ more explicitly (and using the definition of $a$ from \eqref{eqn:a-def}), to prove the second assertion of \eqref{eqn:nn-pre-post-convg}, we want to show
\begin{equation}
\label{eqn:nnac-to-show}
	\frac{\bar{N}_n^{\ac}(k)}{n} \probc \gamma \int_0^a (2+\beta) e^{(2+\beta)u} g_k(a-u) du. 
\end{equation}
For $s\geq 0$, define $\bar{Z}_n^{\ac}(k,s)$ for the number of individuals born in the interval $[\tau_{\gamma n}, \tau_{\gamma n}+s]$ such that by time $\tau_{\gamma n}+s$, these vertices have at least $k$ children. Then note that $\bar{N}_n^{\ac}(k) = \bar{Z}_n^{\ac}(k,\Upsilon_n)$. Mimicking the proof of $N_n^{\bc}(k)$, it is enough to show that 
  \begin{equation}
  \label{eqn:nn-ac-ets}
  	\frac{\bar{Z}_n^{\ac}(k,a_n)}{n} \probc \gamma \int_0^a e^{(2+\beta)u} g_k(a-u) du,
  \end{equation}
where $a_n$ is either the sequence $a_n^{-}= a_n-n^{-1/3}$ or $a_n^{+}=a_n+n^{-1/3}$. To ease notation we will just work with the sequence $a_n=a$. The entire proof goes through by replacing $a$ in the steps below by $a_n$. 

We start with a few preliminary results. The first result describes strong concentration results of the growth of the number of individuals in $\BP_{\mvtheta}^n$ in the interval $[\tau_{\gamma n}, \tau_{\gamma n}+s]$. Define the process 
\begin{equation}
\label{eqn:szn-def}
	\sZ_n(u):= |\BP^{n}_{\mvtheta}(\tau_{\gamma n} +u)|, \qquad 0\leq u\leq a. 
\end{equation}

\begin{prop}
	\label{prop:szn-conc}
	There exists a constant $C<\infty$ such that for all $n$,
	\[\pr\left(\sup_{0\leq u\leq a} |\sZ_n(u) - n\gamma e^{(2+\beta)u}| > \sqrt{n\log{n}}\right) \leq \frac{C}{\log{n}}. \]
\end{prop}
\noindent{\bf Proof:} The plan is to use Doob's $L^2$-maximal inequality for continuous time Martingales (see e.g. \cite[Chapter 1.9]{liptser2012theory}). For this we will need to derive Martingales related to the process $\sZ_n(\cdot)$. Throughout we will write $\set{\cF_t^n: 0\leq t\leq a}$ for the filtration $\set{\BP_{\mvtheta}(\tau_{\gamma n}+t): 0\leq t\leq a}$.  Recall from the rate description in \eqref{eqn:rate-descp} that $\sZ_n(\cdot)$ is a pure birth process such for any $t\geq 0$, conditional on $\cF_t^n$,  $\sZ_n(t)\leadsto \sZ_n(t)+1$ at rate $(2+\beta)\sZ_n(t)-1$. Arguing as in the proof of Lemma \ref{lem:palpha-martingale} it is easy to check that the process 
\begin{equation}
\label{eqn:m1-def}
	M_1(t):= \left(e^{-(2+\beta)t}\sZ_n(t) - n\gamma\right) -\frac{e^{-(2+\beta)t}-1}{2+\beta}, \qquad 0\leq t\leq a,
\end{equation} 
is a mean zero Martingale. This in particular gives that 
\begin{equation}
\label{eqn:szn-mean}
	e^{-(2+\beta)t}\E(\sZ_n(t))= n\gamma + \frac{e^{-(2+\beta)t}-1}{2+\beta}, \qquad 0\leq t\leq a. 
\end{equation}
By Doob's $L^2$-maximal inequality applied to the process $M_1(\cdot)$ we have for any $\lambda > 0$,
\begin{equation}
\label{eqn:doob-m1}
	\pr\left(\sup_{0\leq t\leq a }\left|\left(e^{-(2+\beta)t}\sZ_n(t) - n\gamma\right) -\frac{e^{-(2+\beta)t}-1}{2+\beta}\right|\geq \lambda\right)\leq \frac{\E(M_1^2(a))}{\lambda^2}. 
\end{equation}
If we can show there exists a constant $C< \infty$ such that $\E(M_1^2(a))\leq Cn$, using $\lambda =.5\sqrt{n\log{n}}$ and algebraic manipulation of \eqref{eqn:doob-m1} completes the proof. So let us now derive this bound on $\E(M_1^2(a))$. 

First squaring the expression in \eqref{eqn:m1-def}, expanding and using \eqref{eqn:szn-mean} gives for $t\geq 0$,
\begin{equation}
\label{eqn:szn-square-first}
	\E(M_1^2(t))= \E\left(e^{-(2+\beta)t}\sZ_n(t) - n\gamma\right)^2 -\left(\frac{e^{-(2+\beta)t}-1}{2+\beta}\right)^2. 
\end{equation}
Thus we need to understand the evolution of the process $\sZ_n^2(\cdot)$. Again using the rate description of $\sZ_n$, this process undergoes a change 
\[\Delta \sZ_n^2(t):= \sZ_n^2(t+) - \sZ_n^2(t) = (1+2\sZ_n(t)),\] 
at rate $(2+\beta)\sZ_n(t) - 1$. Using this one may check that the following process on $[0,a]$ 
\begin{equation}
\label{eqn:m2-def}
	M_2(t):= e^{-2(2+\beta)t}\sZ_n^2(t) -  \int_0^t e^{-2(2+\beta)s}\beta \sZ_n(s)ds -   \frac{e^{-2(2+\beta)t}}{2(2+\beta)},
\end{equation}
is also a Martingale. In particular since first moments are conserved, 
\begin{equation}
\label{eqn:szn2-mean}
	\E(e^{-2(2+\beta)}\sZ_n^2(t))= n^2\gamma^2 + \int_0^t \beta e^{-2(2+\beta)s}\E(\sZ_n(s)) ds -\frac{e^{-2(2+\beta)t}-1}{2(2+\beta)}.  
\end{equation}
Using \eqref{eqn:szn-mean} shows that there exists a constant $C$ such that 
\begin{equation}
\label{eqn:szn2-ngamma2}
	\left|\E(e^{-2(2+\beta)}\sZ_n^2(t))- n^2\gamma^2 \right|\leq n\gamma.
\end{equation}
Expanding the first bracket in \eqref{eqn:szn-square-first}, using \eqref{eqn:szn-mean} and \eqref{eqn:szn2-ngamma2} shows that $\E(M_1^2(a))\leq Cn$ for some constant $C$. This completes the proof. 

\qed

Now divide the interval $[\tau_{\gamma n},\tau_{\gamma n}+a]$ into $an^{1/3}$ intervals of length $n^{-1/3}$:
\[\set{\left[\tau_{\gamma n},\tau_{\gamma n}+\frac{1}{n^{1/3}}\right], \left[\tau_{\gamma n}+\frac{1}{n^{1/3}}, \tau_{\gamma n}+\frac{2}{n^{1/3}}\right], \ldots, \left[\tau_{\gamma n}+\frac{an^{1/3}-1}{n^{1/3}},\tau_{\gamma n}+\frac{an^{1/3}}{n^{1/3}}\right]},\]
  of length $n^{-1/3}$.  To ease notation, write the above collection as $\set{\cI_i:0\leq i\leq an^{1/3}-1}$. Further let $\tau_i^n = \tau_{\gamma n} + i/n^{1/3}$ with $\tau_0^n = \tau_{\gamma n}$ so that $\cI_i = [\tau_i^n, \tau_{i+1}^n]$. 
  
  Now write $\birth_i$ for the collection of vertices that were born in interval $\cI_i$ (i.e. the collection of vertices $v$ with birth times $T_v\in \cI_i$) and write \[\sZ_n(\cI_i) := |\birth_i| = \sZ_n\left(\tau_{i+1}^n\right) - \sZ_n\left(\tau_{i+1}^n\right),\]
  for the number of individuals born in this interval. Then the following is an easy corollary of Proposition \ref{prop:szn-conc}. 
  \begin{cor}\label{cor:szni-good-conc}
  	We have 
	\[\pr\left(\bigcap_{i=0}^{an^{1/3}-1}\set{\left|\sZ_n(\cI_i)- (2+\beta)\gamma n^{2/3}e^{\frac{(2+\beta)i}{n^{1/3}}}\right| < 2\sqrt{n\log{n}}}\right)\to 1,\]
	as $n\to\infty$. 
  \end{cor}
  For future reference write $\sG_n$ for the event above namely
  \begin{equation}
  \label{eqn:sgn-def}
  	\sG_n:= \bigcap_{i=0}^{an^{1/3}-1}\set{\left|\sZ_n(\cI_i)- (2+\beta)\gamma n^{2/3}e^{\frac{(2+\beta)i}{n^{1/3}}}\right| < 2\sqrt{n\log{n}}}
  \end{equation}
  
  Now for each interval $\cI_i$,  we will partition the vertices born in this interval into two classes:
 \begin{enumeratea}
 \item The collection of good vertices $\cG_i$: This consists of all $v\in \birth_i$ such that they produce {\bf no} children by the end of the interval i.e. vertices $v$ with $T_v\in [\tau_{\gamma n}+i/n^{1/3},\tau_{\gamma n}+ (i+1)/n^{1/3}]$ such that by time $\tau_{\gamma n}+(i+1)/n^{1/3}$, vertex $v$ still has no children. Note that since the intervals are of time length $n^{-1/3}$, one expects a large proportion of vertices born in the interval $\cI_i$ to be good. Write $\sZ_n^{\sss\good}(\cI_i) = |\cG_i|$ for the number of good vertices in $\cI_i$. 
 \item The collection of bad vertices $\cB_i:= \birth_i\setminus \cG_i$, the collection of vertices born in $\cI_i$ which produce at least one child by time $\tau_{\gamma n}+i/n^{1/3}$. Write $\sZ_n^{\sss\bad}(\cI_i) = |\cB_i|$ for the number of such bad vertices in $\cI_i$.  Write \[\displaystyle \sZ_n^{\sss\bad}:= \sum_{i=0}^{an^{1/3}-1} \sZ_n^{\sss\bad}(\cI_i)\] for the total number of bad vertices.  
 	
 \end{enumeratea} 

Fix a constant $C$ and define the event  $B_i^n = \set{\sZ_n^{\sss\bad}(\cI_i) \geq C n^{1/3}\log{n} }$. These events depend on $C$ but we suppress this in the notation. 

\begin{prop}
	\label{prop:bad-set-bound}
	We can choose constant $C<\infty$ large such that $\displaystyle \pr(\cup_{i=1}^{an^{1/3}} B_i^n)\to 0$ as $n\to\infty$. In particular for the total number of bad vertices we have $\sZ_n^{\sss \bad}= O_P(n^{2/3}\log{n})$. 
	 
\end{prop}

\noindent {\bf Proof:} Fix an interval $\cI_i$. Note that every bad vertex is one of two types:
\begin{enumeratea}
	\item A vertex that is a direct child of a vertex born before this time interval. Write $\cD_n^{\sss\bad}$ for these \emph{direct} bad vertices and write $\sD_n^{\sss \bad}(\cI_i) = |\cD_n^{\sss\bad}|$ for the number of such vertices. Further write $\sD_{n,\star}^{\sss \bad}(\cI_i)$ for the total number of descendants of direct bad vertices born in the interval $\cI_i$ (including the direct bad vertices). 
	\item A vertex that is bad and is a child of a vertex born in $\cI_i$. Thus the parent of this vertex is necessarily bad.  
\end{enumeratea}
Thus in particular we have that $\sZ_n^{\sss \bad}(\cI_i)\leq \sD_{n,\star}^{\sss \bad}(\cI_i)$. Now note that direct bad vertices in $\cD_n^{\sss \bad}$ are created via the following steps:

\begin{enumeratei}
	\item A descendant (maybe good or bad) of a vertex born before $\cI_i$ is born into the system. The number of such individuals $\sR_n(\cI_i)\leq \sZ_n(\cI_i)$, the total number of individuals born in the interval $\cI_i$.  Using Corollary \ref{cor:szni-good-conc}, there exists a constant $C$ such that whp as $n\to\infty$, \emph{for all} the intervals $0\leq i\leq an^{1/3}-1$, $\sR_n(\cI_i) \leq Cn^{2/3}$. 
	\item Conditional on {\bf all these} descendants of vertices born before $\cI_i$, such a descendant has to give birth to one individual in the interval $[i/n^{1/3}, (i+1)/n^{1/3}]$. Recall that the time to give birth to the first child is an exponential random variable $E_1$ with rate $(2+\beta)$.  Thus the probability of birthing this first child is bounded by 
	\[p_n = \pr(E_1 \leq n^{-1/3}) \sim \frac{2+\beta}{n^{1/3}}.  \]
	Further by construction none of these vertices can have a parent child relationship and thus their offspring lineages evolve independently.  
\end{enumeratei}
In particular, conditional on all descendants of vertices born before time interval $\cI_i$, 
\begin{equation}
\label{eqn:sdn-stoch}
	\sD_n^{\sss \bad}(\cI_i)\leq_{\stt}  \text{Bin}(\sR_n(\cI_i), p_n)
\end{equation}
Here $\stt$ denotes stochastic domination. Thus using Corollary \ref{cor:szni-good-conc}, \eqref{eqn:sdn-stoch} and standard tail bounds for the Binomial distribution implies that there exists a constant $C<\infty$ such that  
\begin{equation}
\label{eqn:sdn-bound}
	\pr(\sD_n^{\sss \bad}(\cI_i)\leq C n^{1/3}\log{n}~\forall 0\leq i\leq an^{1/3}-1) \to 1, 
\end{equation}
as $n\to\infty$. 

Let us now complete the analysis of $\sD_{n,\star}^{\sss\bad}(\cI_i)$. Let us start with the evolution of descendants of a single bad \emph{direct} vertex after it gives birth to its child. This process then starts reproducing at rate $2+\beta+ 1+\beta = 3+\beta$. Further whenever a new vertex is added to the system, the rate of production increases by at most $2+\beta$. Thus writing $K = \lfloor 3+\beta \rfloor$ and $\nu = 2+\beta$, the number of descendants of such a bad vertex can be bounded by a rate $\nu$ Yule process (see Definition \ref{def:yule-process}) that starts with $K$ individuals at time zero. Write $\set{Y_{\nu}^{K}(t):t\geq 0}$ for such a process.  Thus the number of descendants of such a bad vertex in the time interval $[\tau_{\gamma n}+ i/n^{1/3},\tau_{\gamma n}+ (i+1)/n^{1/3}]$ can be stochastically bounded by $Y_{\nu}^K(n^{-1/3})$. In particular, conditional on $\cD_n^{\sss\bad}(\cI_i)$, 
\begin{equation}
\label{eqn:cdn-star-st}
	\sD_{n,\star}^{\sss \bad}(\cI_i)\leq_{\stt} \sum_{j=1}^{\sD_n^{\sss\bad}(\cI_i)} Y_{\nu,j}^{K}(n^{-1/3}). 
\end{equation}
Here $\set{Y_{\nu,j}^{K}(\cdot):j\geq 1}$ are an \emph{iid} collection of Yule processes with distribution $Y_{\nu}^K(\cdot)$. Using the explicit distribution of the Yule process at a fixed time (Lemma \ref{lem:yule-prop}), it is easy to check that given constant $C> 0$ we can find $A> 0$ such that 
\begin{equation}
\label{eqn:tail-cdn-star}
	\pr\left(\left.\sD_{n,\star}^{\sss \bad}(\cI_i) \geq 10KCn^{1/3}\log{n}\right|\sD_n^{\sss\bad}(\cI_i)\leq Cn^{1/3}\log{n}\right)\leq \exp(-A n^{1/3}). 
\end{equation}
Using this exponential bound with \eqref{eqn:sdn-bound} completes the proof. 
\qed

We now proceed with the proof of \eqref{eqn:nn-ac-ets}. For $0\leq i\leq an^{1/3}-1$, let $Z_{n,{\sss \good}}^{\ac}(k,a:\cI_i)$ be the number of {\bf good} vertices in $\birth_i$ which have at least $k$ children by time $a$. Then note that conditional on $\BP_{\mvtheta}^n({\tau_{i+1}^n})$, 
\begin{equation}
\label{eqn:zngood-dist}
	Z_{n}^{{\sss \good}}(k,a:\cI_i)\stackrel{d}{=} \text{Bin}\left(\sZ_n^{\sss \good}(\cI_i), g_k\left(a-\frac{i+1}{n^{1/3}}\right)\right). 
\end{equation}
Define the events 
\[G_i^n:= \set{\left|Z_{n}^{{\sss \good}}(k,a:\cI_i) -\gamma(2+\beta) n^{2/3} e^{\frac{(2+\beta)i}{n^{1/3}}} g_k\left(a-\frac{i+1}{n^{1/3}}\right)\right| < C n^{1/3}\log{n} }\]
\begin{prop}
	\label{prop:good-vertex-close}
	There exists a constant $C<\infty$ such that 
	$\pr\left(\cap_{i=1}^{an^{1/3}} G_i^n\right) \to 1$ as $n\to\infty$. 
\end{prop}
\noindent{\bf Proof:} Note that $\sZ_n^{\sss \good}(\cI_i)= \sZ_n(\cI_i) - \sZ_n^{\sss \bad}(\cI_i)$. Combining Corollary \ref{cor:szni-good-conc} with Proposition \ref{prop:bad-set-bound} implies that  
\[\pr\left(\bigcap_{i=0}^{an^{1/3}-1}\set{\left|\sZ_n^{\sss \good}(\cI_i)- (2+\beta)\gamma n^{2/3}e^{\frac{(2+\beta)i}{n^{1/3}}}\right| < 3\sqrt{n\log{n}}}\right)\to 1,\]
Now using the distributional identity \eqref{eqn:zngood-dist} and standard tail bounds for the Binomial distribution completes the proof. 

\qed 

We are finally in a position to complete the proof of \eqref{eqn:nn-ac-ets}. First note that 
\begin{equation}
\label{eqn:zkn-bound}
\displaystyle\sum_{i=0}^{an^{1/3}-1} Z_{n}^{{\sss \good}}(k,a:\cI_i)	\leq \bar{Z}_n^{\ac}(k,a) \leq \sum_{i=0}^{an^{1/3}-1} Z_{n}^{{\sss \good}}(k,a:\cI_i) + \sZ_n^{\sss \bad}.  
\end{equation}
Using Proposition \ref{prop:bad-set-bound} $n^{-1} \sZ_n^{\sss \bad} \probc  0 $. Using Proposition \ref{prop:good-vertex-close}
\begin{align*}
	\frac{\sum_{i=1}^{an^{1/3}} Z_{n}^{{\sss \good}}(k,a:\cI_i)}{n}&\sim \frac{\gamma(2+\beta)}{n^{1/3}} \sum_{i=0}^{an^{1/3}-1} e^{\frac{(2+\beta)i}{n^{1/3}}} g_k\left(a-\frac{i+1}{n^{1/3}}\right)\\
	&\to \gamma(2+\beta) \int_0^{a} e^{(2+\beta)u} g_k(a-u) du. 
\end{align*} 
This completes the proof of \eqref{eqn:nn-pre-post-convg} and thus the assertion of the convergence of the degree distribution of the model to the asserted limit in Theorem \ref{thm:deg-dist}.

\qed

We conclude this Section with a related result regarding the evolution of the degree distribution. This follows by directly modifying the proof above. Recall the definitions of $N_n(k,m)$ and $\hat{N}_n(k,t)$ from Section \ref{sec:res-stat}. For future use define for each $k\geq 1$ and $0\leq t\leq 1$
\begin{equation}
\label{eqn:nn-geq}
	N_{n,\geq}(k,m) = \sum_{j\geq k} N_n(j,m), \qquad \hat{N}_{n,\geq}(k,t) = \sum_{j\geq k} \hat{N}_n(j,t),
\end{equation}
namely the number of vertices with degree at least $k$ respectively at discrete time $m$ and at time $t$ when we rescale time by $n$. Write $\hat{q}_{\geq}^{\sss (n)}(k,t) = \hat{N}_{n,\geq}(k,t)/n $. Note that since we divide by $n$ and not $nt$ in this expression we have $\sum_{k=1}^{\infty} \hat{q}_{\geq}^{\sss (n)}(k,t) = t$.  
Now note that by Lemma \ref{lem:as-pre-convg} we have for each fixed $0< t\leq \gamma$, 
\begin{equation}
\label{eqn:pre-phat-nk-cong}
	\hat{p}^{\sss (n)}(k,t) \probc p_{\alpha}(k)= p^{\sss(\infty)}(k,\gamma),
\end{equation}
where $p_{\alpha}(k)$ as in \eqref{eqn:pk-alpha-def} is the limiting degree distribution with no change point. For $\gamma \leq t\leq 1$, analogous to the definition of $a$ in \eqref{eqn:a-def}  define 
\begin{equation}
\label{eqn:at-def}
	a(t):= \frac{1}{2+\beta}\log\frac{t}{\gamma} 
\end{equation}
Analogous to the definition of $D_{\mvtheta}$ in Section \ref{sec:prelim}, define $D_{\mvtheta}(t)$ by replacing $a$ by $a(t)$ throughout the construction. Thus $D_{\mvtheta} = D_{\mvtheta}(1)$. Let 
\begin{equation}
\label{eqn:pk-after-cp}
	p^{\sss(\infty)}(k,t):= \pr(D_{\mvtheta}(t) = k), \qquad k\geq 1,\;\;\gamma \leq t\leq 1. 
\end{equation}
Let $p^{\sss(\infty)}_{\geq }(k,t) = \pr(D_{\mvtheta}(t)\geq k)$. For $0\leq t\leq 1$, let $q_{\geq}^{\sss(\infty)}(k,t) = t p^{\sss(\infty)}_{\geq }(k,t)$. 
\begin{prop}
	\label{prop:dmvtheta-t-sup}
	 For all $k\geq 1$ we have
	 \[\sup_{0\leq t\leq 1}|\hat{q}_{\geq }^{\sss (n)}(k,t) - q_{\geq}^{\sss(\infty)}(k,t)| \probc 0,\]
	 as $n\to\infty$. 
\end{prop}
\noindent {\bf Proof:} For fixed $t\geq \gamma$, define the stopping time 
\[\tau_{tn} = \inf\set{s: |\BP_{\mvtheta}^n(s)| = tn},\]
namely the first time that the continuous time embedding reaches size $tn$. Note that at this time, the corresponding tree has distribution $\cT_{tn}$. Write $\Upsilon_n(t)=\tau_{tn} -\tau_{\gamma n}$ for the amount of (continuous) time it takes for the process to reach this size after the change point. Then note that by Proposition \ref{prop:szn-conc} we can choose an appropriate constant $C<\infty$ such that 
\begin{equation}
\label{eqn:sup-conc-stop}
	\pr\left(\sup_{\gamma \leq t\leq 1}|\Upsilon(t) - a(t)| \leq C\sqrt{\frac{\log{n}}{n}}\right)\to 1,
\end{equation}
as $n\to\infty$, where $a(t)$ is as defined in \eqref{eqn:at-def}. Repeating the above proof for the convergence of degree distribution and replacing $a$ by $a(t)$ throughout the argument shows that for each $t\geq \gamma$ $\hat{N}_{n,\geq}(k,t)/nt\probc \pr(D_{\mvtheta}(t)\geq k)$. Combining this with \eqref{eqn:pre-phat-nk-cong} implies that we have pointwise convergence $\hat{q}_{\geq}^{\sss(n)}(k,t)\to q_{\geq }^{\sss(\infty)}(k,t)$. Now note that for each fixed $n$, the function $\hat{q}_{\geq}^{\sss(n)}(k,\cdot)$ is non-decreasing on $[0,1]$ while the limit function is also monotonically increasing and continuous (and thus uniformly continuous). Given $\eps>0$, fix $\delta>0$ such that for any $t,s\in [0,1]$ with $|t-s|< \delta$,
\[|q_{\geq }^{\sss(\infty)}(k,t) - q_{\geq }^{\sss(\infty)}(k,s) |< \frac{\eps}{4}. \] 
Divide $[0,1]$ into intervals $\set{[i\delta, (i+1)\delta]}$ for $1\leq i\leq 1/\delta$ of length $\delta$. Via the pointwise convergence above, get $n_0<\infty$ large such that 
for all $n> n_0$
\begin{equation}
\label{eqn:idelta-sup}
	\pr\left(\sup_{1\leq i\leq \frac{1}{\delta}}\left|\hat{q}_{\geq }^{\sss(n)}(k,i\delta) -{q}_{\geq }^{\sss(\infty)}(k,i\delta)  \right|<\frac{\eps}{4}\right)\geq 1-\eps.
\end{equation}
Write $G_n(\eps,\delta)$ for the event in the above equation. Then on this event, by the choice of $\delta$, for all $i$ we have $|\hat{q}_{\geq}^{\sss(n)}(k,i\delta)-\hat{q}_{\geq}^{\sss(n)}(k,(i+1)\delta)|\leq \eps/2$.  Using monotonicity, for any $t\in [i\delta, (i+1)\delta]$, $|\hat{q}_{\geq}^{\sss(n)}(k,i\delta)-\hat{q}_{\geq}^{\sss(n)}(k,t)|\leq \eps/2$. By the triangle inequality on $G_n(\eps, \delta)$, for all $t\in [0,1]$ and $n> n_0$,
\begin{align*}
	|\hat{q}_{\geq}^{\sss(n)}(k,t)-q_{\geq}^{\sss(\infty)}(k,t)|\leq |&\hat{q}_{\geq}^{\sss(n)}(k,t)-q_{\geq}^{\sss(n)}(k,i\delta)|+|\hat{q}_{\geq}^{\sss(n)}(k,i\delta)-q_{\geq}^{\sss(\infty)}(k,i\delta)|\\
	&+|{q}_{\geq}^{\sss(\infty)}(k,i\delta)-q_{\geq}^{\sss(\infty)}(k,t)|\leq \frac{\eps}{2}+\frac{\eps}{4}+\frac{\eps}{4}=\eps. 
\end{align*} 
Since $n_0$ is independent of $t$, this completes the proof. 

\qed

\subsection{Proof of the tail exponent for the limiting degree distribution}
\label{sec:proof-tail-bd}
The aim of this Section is to prove the asserted tail bound, namely \eqref{eqn:mvtheta-tail}. First note that the lower tail bound is obvious since with probability $\gamma$, $D_{\mvtheta}$ stochastically dominates $D_\alpha$ and by \eqref{eqn:tail-bound}, $D_\alpha$ has the asserted tail behavior. The main crux is then proving the upper bound, namely 
\begin{equation}
\label{eqn:up-bound-tail}
	\pr(D_{\mvtheta}\geq x) \leq c^\prime/x^{2+\alpha}.
\end{equation}
Recall Definition \ref{def:yule-process} of the Yule process and in particular Lemma \ref{lem:yule-prop} on finite time marginal distribution of the Yule process.

Now note that in the description of the limit random variable $D_{\mvtheta}$, with probability $1-\gamma$, $D_{\mvtheta}  = N_\beta[0,\age] \leq_{\stt} N_\beta[0,a]$ where as before $\leq_{\stt}$ represents stochastic domination. Now define
\begin{equation}
\label{eqn:nu-K-def}
	\nu=2+\beta, \qquad K=\lfloor 1+\beta \rfloor 
\end{equation}
 Let $Y_{\nu}^K$ be a rate $\nu$ Yule process started with $K$ individuals at time zero. Comparing the rate of production of new individuals in the point process $\cP_\beta$ with $Y_{\nu}^K$,  we get that $N_\beta[0,a] \leq_{\stt} Y_\nu^K(a) $. By Lemma \ref{lem:yule-prop}, $Y_\nu^K(a)$ is the sum of $K$ independent Geometric random variables. Using the fact that a geometric random variable has finite moment generating function in a neighborhood of zero and an elementary Chernoff bound implies that there exist constants $\kappa,\kappa^\prime > 0$ such that for all $x\geq 1$, we have an exponential tail bound,
\begin{equation}
\label{eqn:exp-tail-bound-nb-age}
	\pr(N_\beta[0,\age] > x) \leq \pr(Y_{\nu}^K(a) > x)\leq  \kappa^\prime \exp(-\kappa x),
\end{equation}
Thus when with probability $1-\gamma$ $D_{\mvtheta} = N_{\beta}[0,\age]$ then the corresponding random variable has exponential tail.   Thus the main contribution to the tail arises when with probability $\gamma$, $D_{\mvtheta} = D_\alpha + N_\beta^{D_\alpha}[0,a]$. Arguing as above (and assuming $\beta \geq 1$), conditional on $ D_\alpha = k$, we have 
\[N_\beta^{D_{\alpha}} [0,a] \leq_{\stt} \sum_{j=1}^k Y_\nu^{(j),K}(a), \]
 where $\set{Y_\nu^{(j),K}(\cdot):j\geq 1}$ are a collection of independent rate $\nu$ Yule processes each started at time zero with $K$ individuals  and independent of $D_\alpha$. The following elementary Lemma completes the proof. 
 \begin{lemma}
	\label{lem:sum-exp-heavy-tail}
 	Let $D\geq 1$ be non-negative integer valued random variable with $\pr(D \geq x) \leq c/x^\gamma$ for all $x\geq 1$, for two constants $c,\gamma > 0$. Let $\set{Y_i:i\geq 1}$ be a sequence of independent and identically distributed positive integer valued random variables, independent of $D$. Consider the random variable $D^*:= \sum_{j=1}^D Y_i$.  If $Y_1$ has finite moment generating function in a neighborhood of zero then there exists a constant $c^\prime > 0$ such that for all $x \geq 1$,
	\[\pr(D^* \geq x) \leq c^\prime/x^\gamma. \]
 \end{lemma}
 
 \noindent{\bf Proof:} For the rest of the proof, write $\mu=\E(Y_1)<\infty$. Then note that 
\begin{align*}
\pr(D^* \geq x) &\leq \sum_{j=1}^{\frac{x}{2\mu}} \pr(D=j)\pr(\sum_{i=1}^j Y_i  \geq x) + \pr\left(D \geq \frac{x}{2\mu}\right), \\
&\leq \pr(\sum_{i=1}^{\frac{x}{2\mu}} Y_i  \geq x) + \frac{c}{x^\gamma},	
\end{align*}
where the second equation follows using the fact that $Y_i\geq 1$ for all $i$ and the tail bound for $D$ from the hypothesis of the Lemma. To complete the proof, note that standard large deviation bounds imply (since $Y_i$ has a finite moment generating function about zero) imply that there exist constants $\kappa, \kappa^\prime$ such for all large $x$ 
\[\pr\left(\sum_{i=1}^{\frac{x}{2\mu}} Y_i  \geq x \right) \leq \kappa^\prime \exp(-\kappa x).\]
This completes the proof. 

\qed

The only item left to complete the proof of Theorem \ref{thm:deg-dist} is to show that the change point {\bf does} change the degree distribution from the original (no change point) model. In Section \ref{sec:leafclt-proofs} we will carry out a detailed analysis of the density of leaves which in particular will show that the asymptotic density of leaves $p_{\mvtheta}(1)\neq p_\alpha(1)$.

\subsection{Analysis of the maximal degree}
\label{sec:max-degree-proof}
The aim of this Section is to prove Theorem \ref{thm:max-deg}.
For simplicity and to ease notation, we will deal with $k=1$ namely just the maximal degree. The general case follows in an identical fashion. First note that using \eqref{eqn:max-degree-conv}, writing $M_{\gamma n}(1)$ for the maximal degree of a vertex in $\cT_{\gamma n}$ namely in the tree just before the change point implies that $M_{\gamma n}(1) /n^{1/(2+\alpha)} $ converges weakly to a strictly positive random variable. Since $M_n(1) \geq M_{\gamma n}(1)$, this implies that given any $\eps > 0$, there exists a constant $K^\prime_{\eps} >0$ such that 
\[\liminf_{n\to\infty} \pr\left(\frac{M_n(1)}{n^{1/(2+\alpha)}} > K^\prime_\eps \right) > 1-\eps. \]
Thus to complete the proof of theorem \ref{thm:max-deg} we need to show, given any $\eps> 0$, $\exists K^\prime_\eps < \infty$ such that 
\begin{equation}
\label{eqn:ub-max-degree}
	\limsup_{n\to\infty} \pr\left(\frac{M_n(1)}{n^{1/(2+\alpha)}} < K_\eps \right) \geq 1-\eps.
\end{equation} 
For any vertex $v\in [n]$ time point $m\in [n]$, write $\deg(v,m)$ for the degree of vertex $v$ in $\cT_m$ with the obvious convention that $\deg(v,k) =0$ if $k < v$. Then note that $M_n(1) = \max(M_{\pre}(n),M_{\post}(n))$ where 
\begin{equation}
\label{eqn:pre-post-def}
	M_{\pre}(n):=\max_{v\in [1,n\gamma]}{\deg(v,n)}, \qquad M_{\post}(n):=\max_{v\in [n\gamma+1,n]}{\deg(v,n)}. 
\end{equation}
Let us first analyze the maximal degree of vertices that appeared after the change point. Recall the constant $a$ from \eqref{eqn:a-def} and $\nu, K$ from \eqref{eqn:nu-K-def}. 

\begin{lemma}\label{lem:post-cp-max}
	 We have $\pr(M_{\post}(n) > 2Ke^{\nu(a+1)}\log{n}) \to 0$ as $n\to\infty$. 
\end{lemma}
\noindent{\bf Proof:} We will assume $\beta \geq 1$ below. Else replace $\beta$ with one in the rest of the argument below.  For simplicity write $k_n = 2Ke^{\nu(a+1)}\log{n}$. Recall that in the continuous time embedding, $T_v$ represents the time of birth of vertex $v$ and further for $v\in [\gamma n+1,n]$, each such vertex is equipped with a offspring point process $\cP_\beta^v$. As in Section \ref{sec:proof-tail-bd},  $1+\cP_\beta \leq_{\stt} Y_\nu^K$ where $Y_\nu^K$ is a rate $\nu$ Yule process started with $K$ individuals at time zero. Now note that via our continuous time embedding,
\[M_{\post}(n):= \max_{v\in [\gamma n+1, n]} (1+\cP_\beta^v(0,\tau_n - T_v)),\]   
since by time $\tau_{n}$, a vertex born after the change time has been alive for $\tau_n -T_v \leq \tau_n -\tau_{\gamma n} :=\Upsilon_n$ units of time. Now 
\begin{align}
	\pr(M_{\post}(n) > k_n) &\leq \pr(M_{\post}(n) > k_n , \Upsilon_n < a+1) + \pr(\Upsilon_n > a+1),\notag \\
	&\leq \pr(\max_{v\in [\gamma n+1, n]}(1+ \cP_\beta^v(0,a+1)) > k_n) +\pr(\Upsilon_n > a+1). 
\end{align}
Using Lemma \ref{lem:a-def-asymp} we have $\limsup_{n\to\infty} \pr(\Upsilon_n > a+1) =0$. Let $\set{Y_{\nu,v}^K: v\in [\gamma n+1, n]}$ be a family of independent rate $\nu$ Yule processes started with $K$ individuals at time zero.  Using Lemma \ref{lem:yule-prop} a simple union bound and the choice of $k_n$ implies $\pr(\max_{v\in [\gamma n+1, n]} Y_\beta^v(a+1) > k_n)\to 0$.
\qed

Thus the above Lemma implies that the maximal degree amongst vertices that arrive after the change point is $O_P(\log{n})$. To complete the proof of \eqref{eqn:ub-max-degree}, it is enough to show that \eqref{eqn:ub-max-degree} holds with $M_n(1)$ replaced by $M_{\pre}(1)$. Thus fix $\eps\in (0,1)$. Using Proposition \ref{prop:size-bp-alpha} fix $A = A_\eps$ such that 
\begin{equation}
\label{eqn:a-eps-cond}
	\limsup_{n\to\infty} \pr(\tau_{\gamma n} - \frac{1}{2+\alpha} \log{\gamma n} > A ) \leq \eps/2.
\end{equation}
Now consider the following process $\BP_{\mvtheta,\star}^n$:
\begin{enumeratea}
	\item Run the process $\BP_\alpha$ till time $t_n(A):= \frac{1}{2+\alpha} \log{\gamma n} + A$. 
	\item At this time: all vertices in $\BP_{\alpha}(t_n)$ switch to the dynamics with parameter $\beta$ namely each vertex now reproduces at rate proportional to its out-degree + $1+\beta$. 
	\item Run this process for an additional $a+1$ units of time where $a$ is as in \eqref{eqn:a-def}. 
\end{enumeratea}
Abusing notation, let $M_{\pre,A}^{\star}(1)$ denote the maximal degree by time $t_n +a +1$ of all vertices born {\it before} time $t_n$.  We can obviously couple the original process $\BP_{\mvtheta}^n$ and $\BP_{\mvtheta,\star}^n$ such that on the set $\set{\tau_{\gamma n} - \frac{1}{2+\alpha} \log{\gamma n} \leq A, \Upsilon_n \leq a+1 }$ we have $M_{\pre}(1) \leq M_{\pre,A}^\star(1)$. 

 Further note that for any fixed $K$ we have
\begin{align}
\pr\left(M_{\pre}(1) > K n^{1/(2+\alpha)}\right) \leq &\pr\left(M_{\pre}(1) > K n^{1/(2+\alpha)}, \Upsilon_n < a+1, \tau_{\gamma_n}< \frac{1}{2+\alpha}\log{\gamma n}+A\right) \notag\\
&+ \pr(\Upsilon_n > a+1) + \pr(\tau_{\gamma_n} > \frac{1}{2+\alpha}\log{\gamma n}+A).\notag
\end{align}
First choosing $A$ appropriately as in \eqref{eqn:a-eps-cond} and using Lemma \ref{lem:a-def-asymp} we get that for any fixed $K$, 
\[\limsup_{n\to\infty}\pr(M_{\pre}(1) > K n^{1/(2+\alpha)}) \leq \limsup_{n\to\infty}\pr(M_{\pre,A}^\star(1)> K n^{1/(2+\alpha)}) + \eps/2.\]
The following Lemma completes the proof of \eqref{eqn:ub-max-degree}. 
\begin{lemma}\label{lem:mpre-star-bound}
	Fix $A >0$. Given any $\eps > 0$, we can choose $K= K(A,\eps)< \infty$ such that 
	\[\limsup_{n\to\infty}\pr(M_{\pre,A}^\star(1)> K n^{1/(2+\alpha)}) \leq \eps.\]
\end{lemma}
\noindent {\bf Proof:} First note that till time $t_n(A)$, the process $\BP_{\theta,\star}^n$ is a the continuous time version of a (non-change point) preferential attachment model with attachment parameter $\alpha$. This continuous time embedding was used to derive asymptotics for the maximal degree in \cite{bhamidi2012spectra,bhamidi2012twitter}. In particular the bounds derived in these papers imply the following for a fixed $A$: Write $\tilde{M}_n(1)$ for the maximal degree exactly at time $t_n(A)$. Then there exists $L=L(A,\eps)< \infty$ such that  
\begin{equation}
\label{eqn:tildeMn-bd}
\limsup_{n\to\infty}	\pr(\tilde{M}_n(1) > L n^{1/(2+\alpha)}) \leq \eps/2.
\end{equation}
Now note that on the event $\set{ \tilde{M}_n(1) \leq L n^{1/(2+\alpha)}}$ at time $t_n+a+1$, the degree of every fixed vertex in the system is stochastically dominated by a rate $\nu$ Yule process started with $L n^{1/(2+\alpha)}$ vertices at time zero and run for time $a+1$ where $\nu$ is as in \eqref{eqn:nu-K-def}. Write $D_n$ for such a random variable and note that by the description of the dynamics of the Yule process and Lemma \ref{lem:yule-prop}, we have that 
\begin{equation}
\label{eqn:din-distr}
	D_n \stackrel{d}{=} \sum_{j=1}^{Ln^{1/(2+\alpha)}} {Y_{\nu,j}}(a+1),
\end{equation}
where $\set{Y_{\nu,j}(a+1):j\geq 1}$ are iid Geometric random variables with $p=e^{-\nu(a+1)}$.
 Further note that using Proposition \ref{prop:size-bp-alpha} on the size of the branching process, we can choose $C$ such that 
\begin{equation}
\label{eqn:bptheta-n-bd}
	\limsup_{n\to\infty}\pr(|\BP_{\mvtheta,\star}^n(t_n)| > Cn) \leq \eps/2.
\end{equation}
 Thus on the ``good'' event 
 \[\cG_n:= \set{|\BP_{\mvtheta,\star}^n(t_n)| \leq Cn,\tilde{M}_n(1) \leq L n^{1/(2+\alpha)} },\]
we have that 
\[M_{\pre,A}^\star(1) \leq_{\stt} \max_{1\leq v\leq Cn} D_n^v:= \cM_n\]
where $\set{D_n^v: v\geq 1}$ is an iid sequence with distribution \eqref{eqn:din-distr}. Note that $\E(Y_{\nu,i}(a+1)) = e^{\nu(a+1)}$. Let $K:= 10L e^{\nu(a+1)}$. Then standard large deviations for the Geometric distribution implies that there exists a constant $C^\prime > 0$ such that for all $n\geq 1$
\[\pr(D_n \geq K n^{1/(1+\alpha)}) \leq \exp(-C^\prime n^{1/(1+\alpha)} ).\]  
Thus by the union bound, 
\begin{equation}
\label{eqn:mn-union}
	\pr(\cM_n > K n^{1/(1+\alpha)} ) \leq Cn \exp(-C^\prime n^{1/(1+\alpha)} ) \to 0,
\end{equation}
as $n\to\infty$. 
Thus we have,
\begin{align*}
	\limsup_{n\to\infty}\pr(M_{\pre,A}^\star(1)> K n^{1/(2+\alpha)}) \leq \limsup_{n\to\infty}\pr(\cG_n^c) + \limsup_{n\to\infty} \pr(\cM_n > K n^{1/(2+\alpha)} )\leq \eps,
\end{align*}
using \eqref{eqn:tildeMn-bd}, \eqref{eqn:bptheta-n-bd} and \eqref{eqn:mn-union}. This completes the proof of the Lemma and thus the analysis of the maximal degree asymptotics. 

\qed

\section{Analysis of the proportion of leaves}
\label{sec:leafclt-proofs}

The aim of this Section is to prove Theorem \ref{thm:leaves-fclt}. In the next section we will use the proportion of leaves (degree one vertices) to construct consistent estimators of the change point $\gamma$.  We start in Section \ref{sec:leaf-expec} by deriving strong error bounds between the expected proportion of leaves and the asserted limits in \eqref{eqn:p-inf-zero-form}. Then in Section \ref{sec:leaf-clt-proof} we complete the proof of the functional central limit theorem. We start with some preliminary notation. For the rest of the proof, to ease notation, we will write $N_n(m) := N_n(1,m)$ for the number of leaves in $\cT_{m}$ and let $\hat{N}_n(t) = N_n(nt)$. Recall the asserted limiting proportion $\set{p_t^{\sss(\infty)}: 0\leq t\leq 1}$ from \eqref{eqn:p-inf-zero-form}. For each $n\geq 2$ define the collection of real numbers $\vw_n=\set{w_m:2\leq m\leq n-1}$  
\begin{equation}
\label{eqn:wm-def}
	w_m = \begin{cases}
	\left(1-\frac{1+\alpha}{(2+\alpha)m-1}\right) & \text{ if } 2\leq m\leq n\gamma -1,\\
	\left(1-\frac{1+\beta}{(2+\beta)m-1}\right) & \text{ if } n\gamma \leq m\leq n-1. 
	\end{cases}
\end{equation}
\subsection{Expectation error bounds}
\label{sec:leaf-expec}

The following Proposition is the main result of this Section. 
\begin{prop}\label{prop:explicit-form-p0}
	There exists a constant $C<\infty$ independent of $n$ such that the expectations satisfy 
	\begin{equation}
	\label{eqn:nt-pinf-distance}
		\sup_{n\geq 1} \sup_{0\leq t\leq 1} \left|\E(\hat{N}_n(t)) - nt p^{\sss(\infty)}_t\right|\leq C.
	\end{equation}
	
\end{prop}

%
\begin{rem}
	Note that by Proposition \ref{prop:dmvtheta-t-sup}, we know there exists a function $p^{\sss(\infty)}(0,\cdot)$ such that $\hat{p}^{\sss(n)}(0,t) \to p^{\sss(\infty)}(0,t)$ for $0< t\leq 1$. By the Bounded convergence Theorem,  $\E(\hat{p}^{\sss(n)}(0,t)) \to p^{\sss(\infty)}(0,t) $. Thus the above Proposition implies that $p^{\sss(\infty)}(0,t) = p^{\sss(\infty)}_t$. In particular it shows that the degree distribution owing to the change point is {\bf different} from the degree distribution without change point. This is the final nail in proving Theorem \ref{thm:deg-dist}.  
\end{rem}

\begin{rem}
	 A similar result was shown in the context of no change point in \cite[Section 8.6]{van2009random} and \cite{durrett-rg-book} (not just for leaves but for all fixed $k\geq 1$). Our proof uses slightly different ideas starting from the same point as in \cite{van2009random}. While we do not consider higher degree vertices, as in \cite{van2009random}, the result above can be used as a building block to show identical error bounds for expectations of the number of higher degree vertices about limit constants.  
\end{rem}

\noindent{\bf Proof:}  To ease notation write $\vartheta_n(m) = \E(N_n(m))$. The main crux of the proof is studying a recursion relation for $\vartheta_n(m+1)$ in terms of $\vartheta_n(m)$.  
 We will give a careful analysis of the time period before the change point and then describe how the same ideas give the result for after the change point. 
 
 For each $1< m\leq n$ write $\cL_{m+1}$ for the event that vertex $m+1$ connects to a leaf vertex in $\cT_m$. Then note that conditioning on $\cT_m$, when $m< n\gamma$ we have 
\begin{align}
	\E(N_n(m+1)|\cT_m)&= N_n(m)+1- \pr(\cL_{m+1}|\cT_m)\notag\\
	&=N_n(m)+1-\frac{(1+\alpha)N_n(m)}{(2+\alpha)m-1}\label{eqn:mn-cond-exp}
\end{align}
When $m\geq n\gamma$ we have the same recursion as above but with $\alpha$ replaced by $\beta$. 
Taking full expectations and simplifying gives the following recursion:
\begin{equation}
\label{eqn:mn-recursion}
	N_n(m+1) = 1+w_mN_n(m), \qquad \vartheta_n(m+1)= 1+w_m\vartheta_n(m),
\end{equation}
where $\set{w_m: 2\leq m\leq n}$ are as defined in \eqref{eqn:wm-def}. 

\noindent{\bf Before the change point:}
Repeatedly using this recursion and using the boundary condition $\vartheta_n(2)= 1$ gives for $m+1\leq n\gamma$,
\begin{equation}
\label{eqn:mn-formula}
	\vartheta_n(m+1)= \sum_{s=2}^m \displaystyle\prod_{k=s}^m \left(1-\frac{(1+\alpha)}{(2+\alpha)k-1}\right)
\end{equation}
Now fix $s_0\geq 1$ large enough such that the following three conditions hold:
\begin{enumeratei}
	\item For all $k\geq s_0$
	\[\log{k}+\gamma \leq \sum_{i=1}^k \frac{1}{i} \leq  (\log{k}+\gamma)+\frac{1}{k}.\]
	Here $\gamma$ is the Euler–Mascheroni constant. See \cite{boas1977partial}.
	\item For all $k\geq s_0$, $1-\frac{(1+\alpha)}{(2+\alpha)k-1}\geq 1/2$. 
	\item We may choose a constant $C<\infty$ such that for all $k\geq 1$,
\begin{equation}
\label{eqn:k-1-tail}
	\sum_{i=k}^{\infty} \frac{1}{((2+\alpha)k-1)^2}\leq \frac{C}{k}. 
\end{equation}
	Further there exists a constant $C^{\prime}$ such that for all $s> s_0$,  $|\exp(C/s)-1|\leq C^{\prime}/s$ and \[\left|\left(1-\frac{(1+\alpha)}{(2+\alpha)s -1}\right) - e^{-\frac{(1+\alpha)}{(2+\alpha)s-1}} \right|\leq \frac{C^{\prime}}{s^2}. \]
\end{enumeratei}
 
To ease notation, for the rest of the proof let $\delta= (1+\alpha)/(2+\alpha)$. Using the elementary inequality $1-x\leq e^{-x}$ for $x\in (0,1)$ and the choice of $s_0$ above, the following inequalities with a constant $C=C(s_0,\alpha)<\infty$ are readily verified:
\begin{enumerateA}
	\item For all $m\geq s\geq s_0$,
	\begin{equation}
	\label{eqn:prod-s-m-bound}
		\left|e^{-\sum_{i=s}^m\frac{\delta}{i}} - \left(\frac{s}{m}\right)^{\delta} \right|\leq C\frac{s^{\delta-1}}{m^\delta}. 
	\end{equation}
	\item For all $m\geq s\geq s_0$, 
	\begin{equation}
	\label{eqn:exp-k-1-k}
		\left|e^{-\sum_{i=s}^m\frac{(1+\alpha)}{(2+\alpha)i}} - e^{-\sum_{i=s}^m\frac{(1+\alpha)}{(2+\alpha)i-1}}  \right|\leq C\frac{s^{\delta-1}}{m^\delta}. 
	\end{equation}
	\item For all $m\geq s\geq s_0$,
	\begin{equation}
	\label{eqn:prod-bound}
		\displaystyle\prod_{k=s}^m \left(1-\frac{(1+\alpha)}{(2+\alpha)k-1}\right)\leq C\left(\frac{s}{m}\right)^\delta. 
	\end{equation}
\end{enumerateA}
Now note that by the ``Lindeberg'' trick, for any $s\leq m$ and two collections of non-negative numbers $\set{w_k: s\leq k\leq m}$ and $\set{z_k: s\leq k\leq m}$ we have 
\begin{equation}
\label{eqn:lindeberg}
	\left|\prod_{k=s}^m w_k - \prod_{k=s}^m z_k\right|\leq \sum_{k=s}^m |w_k-z_k| \prod_{s\leq l< k} z_k \prod_{l> k} w_k
\end{equation}
Using this with $w_k= 1-\frac{(1+\alpha)}{(2+\alpha)k-1}$ and $z_k = e^{-\frac{(1+\alpha)}{(2+\alpha)k-1}}$ and using \eqref{eqn:prod-s-m-bound}, \eqref{eqn:exp-k-1-k} and \eqref{eqn:prod-bound} gives the following Lemma. 

\begin{lemma}\label{lem:prd-sm-delta}
	Fix $s_0$ as above.  Writing $\delta =(1+\alpha)/(2+\alpha)$ there exists a constant $C<\infty$ such that for all $m\geq s \geq s_0$, 
	\[\left|\prod_{k=s}^m \left(1-\frac{(1+\alpha)}{(2+\alpha)k-1}\right) - \left(\frac{s}{m}\right)^{\delta}\right| \leq C \frac{s^{\delta-1}}{m^{\delta}}.\]
\end{lemma}
Now using the form of the expectation $\vartheta_n(m)$ in \eqref{eqn:mn-formula}, the error bound in the above Lemma and the integral comparison 
\[\frac{1}{m^{\delta}} \int_{s_0}^{m-1}x^{\delta}dx \leq \sum_{s_0+1}^m \left(\frac{s}{m}\right)^{\delta} \leq \frac{1}{m^{\delta}}\int_{s_0+2}^{m+1} x^{\delta}dx, \]
shows that there exists a constant $C$ such that for $m\leq n\gamma$
\begin{equation}
\label{eqn:vn-bound-prev-cp}
	|\vartheta_n(m) - \frac{m}{\delta}|\leq C.
\end{equation}
This is the assertion for the expected number of leaves before the change point. 

\noindent{\bf After the change point:} We now describe the evolution of $\vartheta_n(m)$ for $n\gamma < m\leq n$. We only give the basic idea as the details are the same as before the change point. First note that by the above analysis, there exists a constant $C$ such that $|\vartheta_n(n\gamma)- n\gamma/\delta|\leq C$. Now the evolution of the process after $\gamma n$ is as in \eqref{eqn:mn-cond-exp} with $\alpha$ replaced by $\beta$. Thus starting at $m> n\gamma$ and using the argument above we get 
\begin{equation}
\label{eqn:mn-recursion-ac}
	\vartheta_n(m+1):= \sum_{s=n\gamma +1}^m \prod_{j=s}^{m} \left(1-\frac{1+\beta}{(2+\beta)j-1}\right) + \vartheta_n(n\gamma)\prod_{j=n\gamma}^m \left(1-\frac{1+\beta}{(2+\beta)j-1}\right)
\end{equation}
Simplifying notation and writing $m=nt$ where $\gamma \leq t \leq 1$ and repeating the arguments above it is easy to check that there exists a constant $C$ independent of $n$ such that 
\begin{equation}
\label{eqn:var-conc-t-after}
	|\vartheta_n(nt) - nt p^{\sss(\infty)}_t| \leq C,
\end{equation} 
where $p^{\sss(\infty)}_t$ is as in \eqref{eqn:p-inf-zero-form}. This completes the proof. 

\qed 

\subsection{Proof of Theorem \ref{thm:leaves-fclt}}
\label{sec:leaf-clt-proof}
A central limit theorem for the number of leaves $N_n(n)$ (in fact all degree counts $N_n(k,n)$) at time $n$ in the setting of no change point was established in \cite{resnick2015asymptotic}.  We will extend this to a functional central limit theorem in the change point setting. First recall the function $\delta_{\alpha}$ from \eqref{eqn:delta-u-def}. 
Define the stochastic process 
\begin{equation}
\label{eqn:mnstar-def}
	M_n^*(t) = \begin{dcases}
		t^{\delta_{\alpha}} \frac{(N_n(nt) -\vartheta_n(nt))} {\sqrt{n}} &\text{ if } t\leq \gamma\\
		\gamma^{\delta_{\alpha}}\left(\frac{t}{\gamma}\right)^{\delta_{\beta}} \frac{(N_n(nt) -\vartheta_n(nt))} {\sqrt{n}} & \text{ if } t\geq \gamma
	\end{dcases}
\end{equation}
Recall the process $M(\cdot)$ in \eqref{eqn:dmt-def} and the relationship between $M$ and $G$. Using Proposition \ref{prop:explicit-form-p0} and the continuous mapping theorem, it is enough to show the following result. 

\begin{prop}
	\label{prop:mn-to-m}
	We have $M_n^*(\cdot) \weakc M(\cdot)$ on $D[0,1]$ as $n\to\infty$.  
\end{prop}
\noindent{\bf Proof:}
The main idea is to study Martingales associated with the $\set{N_n(m): 2\leq m\leq n}$ and then use the Martingale Functional central limit theorem. There are an enormous number of variants of such functional limit theorems under a multitude of conditions. We quote the specific form relevant to this setting. Recall the function $\phi(\cdot)$ and the corresponding diffusion $M(\cdot)$ defined in \eqref{eqn:gt-def}. 
\begin{thm}\label{thm:mart-fclt}[{\cite{durrett1978functional,ethier-kurtz}}]
	For each $n\geq 1$, let $\set{M_n(m): 1\leq m\leq n}$ be a mean zero Martingale with finite second moments adapted to a filtration $\set{\cF_n(m): 1\leq m\leq n}$. Write $\set{X_n(m):1\leq m\leq n}$ for the associated Martingale difference sequence namely $X_n(m) = M_n(m) - M_{n}(m-1)$ with $M_n(0) = 0$.  Assume the following two hypothesis:
	\begin{enumeratei}
		\item For each $0\leq t\leq 1$
		\begin{equation}
		\label{eqn:quad-var-convg}
			V_n(nt):= \sum_{m=1}^{nt} \E\left(\left.\left[X_n(m)\right]^2\right|\cF_n(m-1)\right) \probc \phi(t),\qquad \mbox{as } n\to\infty.
		\end{equation}
		\item For each fixed $\eps> 0$
		\begin{equation}
		\label{eqn:eps-sum-small-mg}
			\sum_{m\leq n} \E\left(\left.\left[X_n(m)\right]^2 \ind\set{|X_n(m)|> \eps}\right|\cF_{n}(m-1)\right)\probc 0. 
		\end{equation}
	\end{enumeratei}
	Then defining the process $\bar{M}_n(t):= M_n(nt)$, one has $\bar{M}_n\weakc M$ in $D[0,1]$. 
\end{thm}
For our example (following \cite{resnick2015asymptotic}) define the process 
\begin{equation}
\label{eqn:nn-star-def}
	N_n^*(m)= \frac{N_n(m)-\vartheta_n(m)}{\prod_{j=2}^{m-1} w_j}, \qquad 2\leq m\leq n.
\end{equation}
Here $w_j$ is as in \eqref{eqn:wm-def}. Using the recursion \eqref{eqn:mn-recursion} results in the following Lemma. 
\begin{lemma}\label{lem:nnstar-mart}
	The process $N_n^*$ is a martingale with respect to the filtration generated by $\set{\cT_m: 2\leq m\leq n}$. 
\end{lemma}
Now define the corresponding Martingale differences $d_n(m) = N_n^*(m) - N_n^*(m-1)$. Define $\Delta_n(m) = \ind\set{m+1\text{ connects to a non-leaf vertex in } \cT_{m-1}}$. Then simple algebra and \eqref{eqn:mn-recursion} implies that for $m\leq n\gamma$
\begin{equation}
\label{eqn:dnm-form}
	d_n(m)= \frac{1}{\prod_{j=2}^{m-1}w_j}\left[\Delta_n(m) + N_n(m-1)\frac{(1+\alpha)N_n(m-1)}{(2+\alpha)(m-1)-1} - 1\right],
\end{equation}
and 
\begin{equation}
\label{eqn:e-deltnm}
	\E(\Delta_n(m)|\cT_{m-1})= 1- \frac{(1+\alpha)N_n(m-1)}{(2+\alpha)(m-1)-1}
\end{equation}
For $m\geq n\gamma$ we have identical formulae as \eqref{eqn:dnm-form} and \eqref{eqn:e-deltnm} but now $\alpha$ is replaced by $\beta$. For the rest of the argument we will replace the denominator for the second term $(2+\alpha)(m-1)-1$ by $(2+\alpha)(m-1)-1$. It is easy to check that the error is negligible and will ease presentation. 

Now use Proposition \ref{prop:dmvtheta-t-sup} which allows us to uniformly approximate $N_n(m-1)/(m-1)$ by $p^{\sss(\infty)}_{m/n}$. Further the asymptotics of $\prod_{j=2}^m w_j$ derived in the previous Section implies that for $m\leq n\gamma$, $\prod_{j=2}^m w_j\sim m^{-\delta_\alpha}$ while for $m> n\gamma$, $\prod_{j=2}^m w_j\sim (n\gamma)^{-\delta_\alpha}(m/n\gamma)^{-\delta_\beta}$ where $\delta_\alpha, \delta_\beta$ as as defined in \eqref{eqn:delta-u-def}. Taking conditional expectations in \eqref{eqn:dnm-form}, using \eqref{eqn:e-deltnm} and using the above approximations results in   
\begin{equation}
\label{eqn:dnm-sq-expec}
	\E([d_n(m)]^2|\cT_{m-1})\sim \begin{cases}
		m^{2\delta_\alpha}\left[\delta_{\alpha}p_{m/n}^{\sss(\infty)}(1-\delta_{\alpha}p_{m/n}^{\sss(\infty)})\right] & \text{ if } m\leq n\gamma, \\
		& \\
	(n\gamma)^{2\delta_\alpha}\left(\frac{m}{n}\right)^{2\delta_{\beta}} \left[\delta_{\beta}p_{m/n}^{\sss(\infty)}(1-\delta_{\beta}p_{m/n}^{\sss(\infty)})\right] & \text{ if } m\geq n\gamma
	\end{cases}
\end{equation}
Now consider the Martingale 
\begin{equation}
\label{eqn:mart-mnm-def}
	M_n(m):= \frac{1}{n^{\delta_\alpha+1/2}} \frac{N_n(m)-\vartheta_n(m)}{\prod_{j=2}^{m-1} w_j}, \qquad 2\leq m\leq n.
\end{equation}
We will apply Theorem \ref{thm:mart-fclt} to this Martingale. Let $\set{X_n(m):2\leq m\leq n}$ denote the corresponding Martingale differences. First fix $t\leq \gamma$ and recall the definition of the cumulative conditional variance $V_n(nt)$ till time $t$ in \eqref{eqn:quad-var-convg}. Using the first expression in  \eqref{eqn:dnm-sq-expec} we get 
\begin{align*}
	V_{n}(nt)&\sim \frac{1}{n^{2\delta_\alpha+1}} \sum_{j=1}^{nt} j^{2\delta_\alpha} \left[\delta_{\alpha}p_{m/n}^{\sss(\infty)}(1-\delta_{\alpha}p_{m/n}^{\sss(\infty)})\right]\\
	&\to \int_0^t s^{2\delta_{\alpha}}\left[\delta_{\alpha}p_{s}^{\sss(\infty)}(1-\delta_{\alpha}p_{s}^{\sss(\infty)})\right]ds = \phi(t),
\end{align*} 
as $n\to\infty$. Thus \eqref{eqn:quad-var-convg} is satisfied for $t\leq \gamma$. A similar calculation now incorporating the second expression in \eqref{eqn:dnm-sq-expec} implies that \eqref{eqn:quad-var-convg} is satisfied for all $t\in [0,1]$ with $\phi$ as in in \eqref{eqn:gt-def}.  Now let us check the second condition namely \eqref{eqn:eps-sum-small-mg}. Note that for $m\leq n\gamma$, $X_n(m)\geq \eps$ implies that $3m^{\delta_{\alpha}}\geq \eps n^{\delta_{\alpha}+1/2}$. For large $n$ this is impossible for all $m\leq n\gamma$. A similar calculation for $m> n\gamma$ completes the proof of \eqref{eqn:eps-sum-small-mg}. Using Theorem \ref{thm:mart-fclt} we get that $M_n(n\cdot)\weakc M(\cdot)$ in $D[0,1]$. Using the asymptotics for $\prod_{j=2}^m w_j$ derived in Section \ref{sec:leaf-expec} in \eqref{eqn:mart-mnm-def} now completes the proof of Proposition \ref{prop:mn-to-m} and thus Theorem \ref{thm:leaves-fclt}. \qed

\section{Consistency of the estimator}
\label{sec:consistency}
The aim of this Section is to prove Theorem \ref{thm:cp-estimator-cons}. Fix a truncation level $\eps> 0$ from zero as in the Theorem.  Recall the time averaged proportion of leaves before and after each time $t$ namely \eqref{eqn:ac-def} and \eqref{eqn:bc-def}. Also recall the expression for the limiting proportion of leaves from \eqref{eqn:p-inf-zero-form}. For any fixed interval $[s,t]\subseteq [0,1]$, define $H[s,t]$ by
\begin{equation}
\label{eqn:gk-def}
	H[s,t]:= \frac{1}{t-s}\int_s^t p_u^{\sss(\infty)} du. 
\end{equation}  
The interpretation is as follows: the above gives the expected proportion of leaves in the large network limit if one were to sample a time point $U\in [s,t]$ uniformly at random. Now define the two functions $\bcpmf^{\sss(\infty)}$ and $\acpmf^{\sss(\infty)}$ via:  
\begin{enumeratea}
	\item {\bf Case 1:} For $\eps\leq t\leq \gamma$ 
	\[\bcpmf^{\sss(\infty)} := p_{\gamma}^{\sss(\infty)},\qquad \acpmf^{\sss(\infty)} := \frac{\gamma-t}{1-t} p_{\gamma}^{\sss(\infty)} + \frac{1-\gamma}{1-t} H[\gamma,1]\]
	\item {\bf Case 2:} For $t> \gamma$
	\[\bcpmf:= \frac{\gamma-\eps}{t-\eps} p_{\gamma}+ \frac{t-\gamma}{t-\eps} H([\gamma,t]), \qquad \acpmf:= H([t,1]).  \]
\end{enumeratea}
Define analogously to \eqref{eqn:dnt-fn} the function
\begin{equation}
\label{eqn:dnt-fn-inf}
	D(t):= (1-t)|\bcpmf^{\sss(\infty)} - \acpmf^{\sss(\infty)}|, \qquad t\in [\eps,1]. 
\end{equation} 
Routine algebra shows that 
	\begin{equation}
	\label{eqn:dt-expression}
		D(t):=\left\{\begin{array}{lc}
			(1-\gamma) \big|p_{\gamma}^{\sss(\infty)}- H[\gamma,1]\big| & \mbox{ for } \eps\leq t\leq \gamma.\\
			& \\
			(1-\eps)\big|H[\eps,t]-H[\eps,1]\big| & \mbox{ for } t > \gamma.
		\end{array} \right.
	\end{equation}
Using the form of the limit proportion $p_t^{\sss(\infty)}$ from \eqref{eqn:p-inf-zero-form} the following result is easy to check. 
\begin{lemma}\label{lem:d-inf-properties}
	Fix $\eps <\gamma$ and assume $\alpha\neq \beta$. Then $D(\cdot)$ is a continuous function on $[\eps,1]$ such that $D(\cdot)$ is constant on the interval $[\eps,\gamma]$ and then is strictly monotonically decreasing on the interval $[\gamma,1]$ with $D(t)\to 0$ as $t\to 1$. Further the function has a strictly negative right derivative at $\gamma$ namely 
	\begin{equation}
	\label{eqn:right-deriv}
		\partial_+ D(\gamma):= \lim_{t\downarrow \gamma}\frac{D(t) -D(\gamma)}{t-\gamma} < 0.
	\end{equation} 
\end{lemma}

Now Theorem \ref{thm:leaves-fclt} immediately results in the following result. 
\begin{lemma}\label{lem:dinf-dn}
	Fix $\eps > 0$. Then 
	\[\sup_{t\in [\eps,1]}|D_n(t) - D(t)| = O_P\left(\frac{1}{\sqrt{n}}\right).\]
\end{lemma}

Now combining Lemmas \ref{lem:d-inf-properties} and \ref{lem:dinf-dn} completes the proof. \qed

\section*{Acknowledgements}
SB has been partially supported by NSF-DMS grants 1105581 and 1310002 and SES grant 1357622. ABN has been partially supported by NSF-DMS grant 1310002 and NIH R01 MH101819-01. JJ has been partially supported by NSF-DMS grant 1310002. We thank Edward Carlstein for several illuminating conversations on change point detection. We thank Frances Tong for extensive help with regards to visualization of various plots in the paper.



\end{document}